\documentclass{article}
\usepackage[margin=1in]{geometry}
\usepackage{amsmath,amsthm,amssymb}
\usepackage{graphicx}
\usepackage{color}
\usepackage[dvipsnames]{xcolor}
\usepackage{subcaption}
\usepackage[export]{adjustbox}
\usepackage{algorithm}
\usepackage{slashbox}
\usepackage[noend]{algpseudocode}
\usepackage{url}

\newcommand{\norm}[1]{\left\|#1\right\|}

\newcommand{\ytrev}[1]{\textcolor{black}{#1}}

\graphicspath{ {./pic/} }

\DeclareMathOperator*{\argmin}{\text{argmin}}

\usepackage{lineno}

\begin{document}


\title{Efficient Deep Learning Techniques for Multiphase Flow Simulation in Heterogeneous Porous Media}

\author{
Yating Wang \thanks{Department of Mathematics, Purdue University, West Lafayette, IN 47907, USA (\texttt{wang4190@purdue.edu})} \and
Guang Lin \thanks{Department of Mathematics, Purdue University, West Lafayette, IN 47907, USA \& Department of Mechanical Engineering, Purdue University, West Lafayette, 47907, USA (\texttt{guanglin@purdue.edu}) } \footnote{Corresponding author.}
}

\maketitle

\begin{abstract}
We present efficient deep learning techniques for approximating flow and transport equations for both single phase and two-phase flow problems. The proposed methods take advantages of the sparsity structures in the underlying discrete systems and can be served as efficient alternatives to the system solvers at the full order. In particular, for the flow problem, we design a network with convolutional and locally connected layers to perform model reductions. Moreover, we employ a custom loss function to impose local mass conservation constraints. This helps to preserve the physical property of velocity solution which we are interested in learning. For the saturation problem, we propose a residual type of network to approximate the dynamics. Our main contribution here is the design of custom sparsely connected layers which take into account the inherent sparse interaction between the input and output. After training, the approximated feed-forward map can be applied iteratively to predict solutions in the long range.  Our trained networks, especially in two-phase flow where the maps are nonlinear, show their great potential in accurately approximating the underlying physical system and improvement in computational efficiency. Some numerical experiments are performed and discussed to demonstrate the performance of our proposed techniques.
\end{abstract}


\section{Introduction}

Various physical phenomena in engineering applications are described by flow and transport problem in porous media, including reservoir engineering, climate dynamics, material science and so on. The underlying problems naturally exhibit heterogeneities covering from large physical scales down to micro-scales. Numerical simulations for these problems are challenging due to a rich class of length scales and potential uncertainties. In the past decades, numerous model reduction techniques and multiscale methods are proposed to design alternative models with great computational efficiency as well as desired accuracy. Local model reduction techniques typically involve building representations of underlying heterogeneity using some basis functions or effective coarse grid properties, and then constructing a form of the coarse level
equations \cite{ab05, egw10,  arbogast02, GMsFEM13, AdaptiveGMsFEM, brown2014multiscale, ElasticGMsFEM, ee03, abdul_yun, fish2004space, fish2008mathematical, oz07, matache2002two, henning2009heterogeneous, OnlineStokes, chung2017DGstokes,WaveGMsFEM, MsDG, NLMC}. For large scale dynamical system, global reduced order models adopting Proper Orthogonal Decomposition method, Krylov subspace projection method, etc, are proposed to approximate the state-space problems. Though both local and global model reduction techniques have been extensively applied in many problems, reduced order models may have complicated forms in the linear case, not mentioning in nonlinear settings \cite{ehg04, nonlinear_AM2015, cegg13, yang2016PodDeim}.

Recently, deep learning has attracted growing attention in a rich class
of applications. It has gained revolutionary results in image, speech, and text recognition \cite{dcnn_im2012, dnn_speech_2012, HK_resnet2016}. The potential of deep neural networks lies in their great capacity in approximating high-dimensional nonlinear maps. There are extensive efforts devoted to learning the expressivity
of deep neural nets theoretically, just to mention a few \cite{Cybenko1989, Hornik1991, Csaji2001, Telgrasky2016, Poggio2016}, where the universal approximation property are investigated and the approximation ability of deep networks to a rich classes of functions are shown. This motivates a number of works bringing deep neural networks to the field of scientific computing. Some recent works include numerically solving parametric partial differential equations \cite{Ying_paraPDE}, designing multigrid solvers \cite{ Fan2018MNNH}, constructing efficient multiscale models \cite{wang2018Gmsfem}, learning surrogate models by deep convolution networks \cite{deepconv_Zabaras, E_deepRitz}, incorporating reduced order models in conjunction with observation data \cite{wang2018NLMC_DNN}, solving ODE systems driven by data \cite{NeuralODE, QinXiu2018dataODE}, \ytrev{combining deep learning concepts and some proper orthogonal decomposition (POD) model reduction methods \cite{pod_dl_SW, pod_dl_zhang}} and so on. In particular, the authors in \cite{wang2018NLMC_DNN} propose novel methods to learn the nonlinear feed-forward map considering the sparsity nature of the underlying problem, and authors in \cite{Fan2018MNNH} also deals with nonlinear maps issued from underlying physical problems by designing multiscale neural networks.

The focus of this work is to apply sparse learning ideas on the coupled flow and transport systems in both single phase and two-phase settings. Our contribution is to design appropriate deep neural networks as computational efficient alternative models to the complicated physical problems.  We start with the single-phase case, where the flow problem is described by Darcy's law and the transport dynamics is driven by a velocity field. Motivated by ideas from model reduction algorithms, we design a neural network architecture to learn the maps between some input (such as the source term or the permeability field) and output (velocity solutions) in the flow problem. In this flow approximation task, we will first transform the input to some coarse features by using convolutional layers or average pooling layers. It \ytrev{is} worth mentioning that, although traditional multi-layer perceptron models with densely connected layers were successfully utilized in many applications, they suffer from the curse of dimensionality because of the full connectivity between layers. To handle this issue, it's natural to come up with some locally connected networks which can benefit from the fact that neighboring nodes have an inherent relationship. For instance, convolutional layers/ locally connected layers have gained remarkable success and performed superior to others. Hence, in our network design for flow approximation, some locally connected layers are adopted. They \ytrev{do} not only take into account the spatial structure of data but also provide sufficient learning capacity to extract rich hidden features. We remark that the locally connected layer is similar to the convolutional layer, but a different set of filters is applied at each different patch of the input. After some coarse features are extracted from the input layer, we then consider working on the coarse level deeply to learn more hidden properties, which is inspired by the analogy to the multiscale method. In the end, a decoding step is realized by a single densely connected layer between the coarse features and the fine velocity solution. Another important aspect of our proposed network is that it preserves the physical property of the approximating quantity. The velocity field used in the training process is locally mass conservative. One expects the predicted solution produced from the trained neural network can also maintain the desired physical property. We tackle this issue by imposing constraints with physical meaning in the loss function while training the network. Our numerical results show that the proposed framework did help to improve the local mass conservative property among predictions. As for the transport equation, we are interested in learning the dynamics between saturation solutions at different time steps. To be specific, the solution state at the time $n+1$ depends on the solution at the time step $n$ and other parameters, such as velocity fields and source terms. It is natural to employ a residual network structure in approximating the feed-forward map. In this approximation task, we still want to rely on the sparsity of the discrete system to define unknowns and design sparse custom network layers. More specifically, our contribution in this paper is to design the Sparse Velocity Layer, which treats the velocity input as multipliers of the weight matrices and also take advantages of the data sparsity structure of the discrete system to create sparse weight matrices. This will result in a notable reduction in the number of trainable parameters for the weight matrices and require less training data due to the simplification of the network. We present detailed network architecture as well as several numerical experiments. The numerical results illustrate the capability of the proposed network and implies both the accuracy and efficiency of our methods.

Solving two-phase system is challenging due to its nonlinearity. We apply our experiences in designing neural networks for single-phase flow and transport approximations, and appropriately generalize them to the two-phase problems. In classical numerical methods, the two-phase flow and saturation equations are solved sequentially. Since the flow velocity not only depends on the absolute permeability but also the relative permeability, which is a function of the phase saturation. Thus one needs to update the velocity field at every time step. Furthermore, the saturation equation is also nonlinear and some iterative methods are required at each time step to obtain the solution. With this regard, we expect to alleviate the computational burden using deep neural networks. The ideas of constructing approximation maps for velocity problem and saturation dynamics are based on the previously proposed architectures. For flow approximation, our main contribution is to develop the input-output locally connected network between the total mobility (a function of saturation) and the velocity solution. For the nonlinear saturation approximation, one of the difficulties is that the velocity field varies in every saturation forward map. However, it's not trivial to utilize every velocity input as multipliers for the weight matrices, the same ways as we did in the single flow case. A workaround is that we compute the mean of the velocity fields in each training batch, and feed in the mean velocity field during the training process of the nonlinear saturation map. The sparse structure of the system is still realized when designing the weight matrices. After both networks for approximating flow and saturation equations are trained, we propose a sequential algorithm to recursively employ the trained networks to predict saturation solutions in the long range. Our numerical results show great potential both in accuracy and computational improvement for approximating the two-phase flow system.

The paper is organized as follows. In Section \ref{sec:single_phase}, we present deep neural network construction in the single-phase flow, as well as corresponding numerical experiments. In Section \ref{sec:two_phase}, the generalization, and adjustment of the proposed methods to two-phase flow problems are illustrated, numerical validations are followed after the methodologies. A conclusion is presented in Section \ref{sec:conclusion}.

\section{Single phase flow}\label{sec:single_phase}

In the single phase flow, we consider \ytrev{solving the flow problem by the mixed finite element method (MFEM). One can take $\text{RT}_0$ (the lowest order Raviart-Thomas element) to approximate velocity, and $P_0$ (piecewise constant element) to approximate pressure. It is well known that MFEM is locally mass conservative. More precisely,} we aim to solve the velocity $u$ and pressure $p$ from a mixed system
\begin{align*}
\kappa^{-1} u + \nabla p &= 0 \quad \quad \text{in}  \quad D\\
\text{div} (u) &= f \quad \quad \text{in}  \quad D\\
u\cdot n &= 0 \quad \quad \text{on}  \quad \partial D
\end{align*}
where $\kappa$ denotes the heterogeneous permeability field (see Figure \ref{fig:perm} as an example), $D$ represents the computational domain, and $f$ is source term.

The fine grid of the problem reads
\begin{align*}
a(u,v)+b(v,p) &= 0   & \text{ for all }  v\in V_h\\
b(u,q) &= -(f,q)  &\text{ for all }  q\in Q_h
\end{align*}
where $a(u,v) = \int_{D} \kappa^{-1} u \cdot v$, and $b(v,p) = -\int_D p \; \text{div} v $. $V_h$ and $Q_h$ are the finite element space for velocity and pressure, respectively.

In addition to the flow equation, the transport equation of saturation $S$ is given by
\begin{equation*}
\displaystyle{ \frac{\partial S} {\partial t} + u \cdot \nabla S = r}
\end{equation*}
where $u$ is the velocity field obtained in the flow problem.

The transport of saturation can be solved by finite volume method on the fine grid. One can discretize the time derivative using the forward Euler scheme, where the time step size can be chosen from the CFL condition. To be specific, on a fine grid $K_i$, the value $S_i$ at time $t^{n+1}$ can be obtained by
\begin{equation} \label{eq:sat_linear}
|K_i| \displaystyle{ \frac{S_i^{n+1} - S_i^{n}}{\text{d} t} + \sum_{e_j \in \partial K_i} F_{ij} (S^n)  = |K_i| r_i}
\end{equation}
where $F_{ij}$ is the upwind flux, i.e.
$$
F_{ij}(S^n) = \left\{
                \begin{array}{ll}
                  \int_{e_j} (u_{ij} \cdot n) S_i^n\quad \text{ if }\quad u_{ij} \cdot n  \geq  0\\
                  \int_{e_j} (u_{ij} \cdot n) S_j^n \quad\text{ if } \quad u_{ij} \cdot n  <  0
                \end{array}
              \right.
$$
Here $e_j$ denotes the edge shared by fine grids $K_i$ and $K_j$, $u_{ij}$ the velocity on the edge $e_j$.

The above flow and transport problem is solved sequentially. One first solves the flow equation to obtain the velocity, and the obtained velocity filed is used to drive the saturation. 
However, once the source term $f$ or the permeability changed, the system needs to be solved again, which is computationally expensive. We are interested in developing efficient alternatives for the high fidelity model and obtain the desired accuracy with reduced computational effort.

\subsection{Model reduction for flow problem using neural network} \label{sec:vel_single} \label{sec:vel_lnn}

\subsubsection{Discrete flow problem solver}

We note that the discrete flow problem has the following matrix formulation on the fine grid $\mathcal{T}_h$
 \begin{equation}\label{eq:vel_mat_single}
 \begin{bmatrix}
    A_h(\kappa) & B_h^T   \\
    B_h & 0
  \end{bmatrix}
   \begin{bmatrix}
    u_h \\
   p_h
  \end{bmatrix} =
   \begin{bmatrix}
   0   \\
  -F
  \end{bmatrix}
 \end{equation}

In practice, one may need to solve the system with thousands of different $\kappa$ (when there are uncertainties in the permeability field) or many different source terms $f$. Typically, the linear system \eqref{eq:vel_mat_single} generated from MFEM on the fine grid can be solved using a direct solver, however, \ytrev{the computational expense grows heavily with the complexity of heterogeneity in the media and number of simulations.}

To tackle this difficulty, one way is to solve the problem on a coarser grid $\mathcal{T}_H$. There are many mixed multiscale methods in literature, for example \cite{MixedGMsFEM, Arbogast_mixed_MS_11, ae07}. We assume the multiscale basis are computed for velocity in each local coarse region, and the piecewise constant on the coarse blocks are used for pressure. Let $R_u$ be the matrix with size $N_u^H \times N_u^h$ consisting of multiscale velocity basis in each row, where $N_u^H$ and $N_u^h$ are the coarse and fine degree of freedom (DOF) for velocity, respectively. Let $R_p$ be the matrix with size $N_p^H \times N_p^h$, which works as an averaging of fine-scale pressure unknowns over each coarse block. $R_u$ and $R_p$ together form a restriction matrix. And the coarse scale system has the form
 \begin{equation}\label{eq:vel_mat_single_coarse}
 \begin{bmatrix}
    A_H & B_H^T   \\
    B_H & 0
  \end{bmatrix}
   \begin{bmatrix}
    u_H \\
   p_H
  \end{bmatrix} =
  \begin{bmatrix}
    R_u & 0   \\
    0 & R_p
  \end{bmatrix}
   \begin{bmatrix}
    A_h(\kappa) & B_h^T   \\
    B_h & 0
  \end{bmatrix}
    \begin{bmatrix}
    R_u^T & 0   \\
    0 & R_p^T
  \end{bmatrix}
   \begin{bmatrix}
    u_H \\
   p_H
  \end{bmatrix}=
   \begin{bmatrix}
   0   \\
  -F_H
  \end{bmatrix}
 \end{equation}
One then can solve the coarse scale system \eqref{eq:vel_mat_single_coarse}, for instance, using the Arrow-Hurwicz iterative method
\begin{equation} \label{eq:iterative_velc}
 \begin{aligned}
u_H^{k+1} &= u_H^{k}+\alpha(A_H u_H^{k}  -B_H^T p_H^k) \\
p_H^{k+1}  &= p_H^k + \omega (B_H u_H^{k+1} + F_H)
\end{aligned}
\end{equation}
where $\alpha > 0$, $\omega > 0$ are relaxation parameters.

Let the initial guesses $p_H^0$ and $u_H^0$ be zero, then after one iteration, we observe that $p_H$ can be obtained by multiplying $F_H$ by some matrix. Next, we obtain $u_H$ by plugging in $p_H$ in the first equation of \eqref{eq:iterative_velc}. This indicates we can view $p_H$ as some inherent unknowns and iteratively obtain the velocity solution on the coarse grid. Next, one can obtain the fine scale velocity by a downscaling step $u_h = R_u^T u_H$.

We remark that $p_H$ will be some hidden coarse grid properties in our deep learning algorithm. In the end, we are interested in the velocity solution only.

\subsubsection{Network architecture}
The focus of this work is to utilize a deep neural network on the approximation of problems of interest. Before that, let's first briefly introduce some basics of deep learning.

Generally, in deep learning, let the function $\mathcal{N}$ be a network of $d$ layers, the $i$-th layer is denoted by $l_i$ ($i = 1, \cdots, m$). Let $x$ be the input and $y$ be the corresponding output.
We write
\begin{equation*}
\mathcal{N}(x; \theta) = \sigma(l_d \sigma (\cdots  \sigma(l_2  \sigma(l_1( x) ) \cdots   ) )
\end{equation*}
where $\sigma$ is the activation function. Denote by $\theta$ all the trainable parameters in the network.
Suppose we are given a collection of sample pairs $(x_j, y_j)$.
The goal is to find $\theta^*$ by solving an optimization problem
\begin{equation*}
\theta^* = \argmin_{\theta} \frac{1}{N}\sum_{j=1}^{N} \mathcal{L}(y_j, \mathcal{N}(x_j; \theta)),
\end{equation*}
where $N$ is the number of the samples. This implies that one aims to find parameter $\theta^*$ which can minimize the mean loss $\mathcal{L}(y_j, \mathcal{N}(x_j; \theta))$ among all training samples using stochastic gradient descent iteratively. Then, the trained network $\mathcal{N}(x; \theta^*) $ will be applied to make predictions on new inputs $x_{\text{new}}$.

To improve the efficiency of solving the flow problem, we are interested in learning the discrete velocity solution $u_h$ (outputs) given different permeability fields $\kappa$ or sources $f$ (inputs). From now on, we omit the subscript $h$ and just use the simple notation $u$ instead. That is, we would like to use pre-computed samples pairs to train a neural network which approximates the map $u = \mathcal{N}(\kappa; \theta)$, or $u = \mathcal{N}(f; \theta)$. When given a test case, one can use the trained network to predict the quantity of interest effectively.

In the single phase flow case, we consider the map between $f$ and $u$ which is linear, a brute-force way is to use one dense layer to connect the input and output. However, as mentioned before, due to the multiple scales of the underlying permeability (see Figure \ref{fig:perm} for illustration), one needs to use the sufficiently fine grid in order to resolve all scale property, which results in large dimensions of $f$ and velocity solution on the fine grid. Then a direct dense connection between the input and output will not be affordable.

With this background in mind, we will propose a network architecture which is deeper and aims to perform the model reduction as described in the previous section. The proposed network will have a small number of trainable parameters compared with a fully connected neural network.

One of the core ingredients we will adopt in the network architecture is the locally connected layer. Let $\alpha$ be the number of filters, each filter has a size $k \times k$, $s$ is the stride size, then the connection between two tensors \ytrev{$Y^1$ and $Y^2$} by a locally connected layer as follows
\ytrev{
\begin{equation}
\displaystyle{ Y^1_{i;k} = \sigma( \sum_{j=(i-1)s+1}^{(i-1)s+k^2} \sum_{l=1}^{\alpha} W_{i,j;k,l}  Y^2_{j;l}   + b_{i;k} ) }
\end{equation}
}
The locally connected layer shares some similarities with the convolutional layer, the difference lies in that a different set of filters is applied at each different patch of the input.

\begin{figure}[!hbt]
\centering
\includegraphics[scale=0.35]{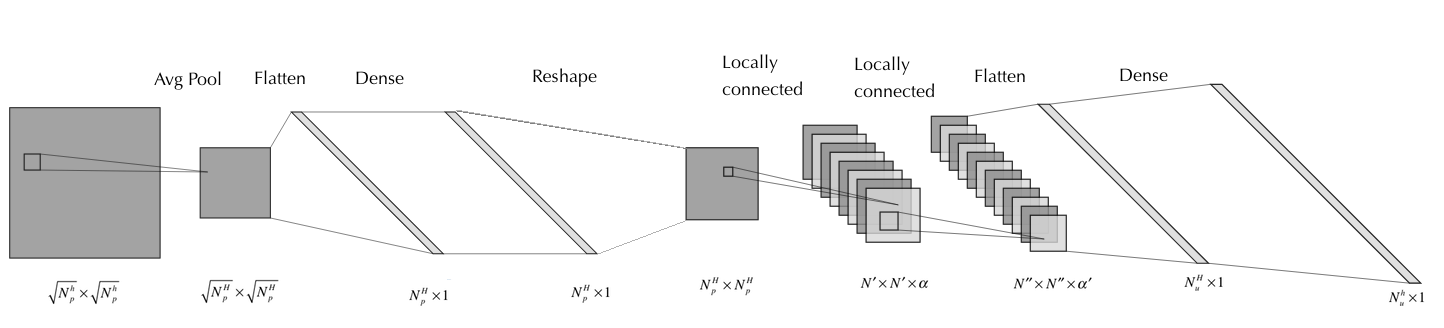}
\caption{An illustration of the network architecture for flow approximation.}
\label{fig:vel_lcn}
\end{figure}

The proposed network consists of the following components. One first project the input tensor $f$ on a coarser grid using an average pooling layer. For example, the input dimension of $f$ is $\sqrt{N_p^h} \times \sqrt{N_p^h}$, where $1/\sqrt{N_p^h}$ is the size of fine grid in the computational domain. We then let the size of a coarse block to be $1/\sqrt{N_p^H}$. Then the average pooling layer are specified to have a pooling size of $\big(\sqrt{N_p^h}/ \sqrt{N_p^H}\big) \times \big(\sqrt{N_p^h}/ \sqrt{N_p^H} \big)$, which will transform the input from the fine grid level to the coarse grid level.  Next, a flattened layer followed by a dense layer with $N_p^H$ neurons is applied to the intermediate output from the previous step. In the dense layer, a square weight matrix is multiplied to its input, which results in some hidden coarse grid property with the same dimension $N_p^H \times 1$. Further, after reshaping the hidden coarse grid property to $\sqrt{N_p^H} \times \sqrt{ N_p^H}$, we use a few locally connected layers to dig in the hidden features. Then, we flatten the resulting hidden features and connect them with the neurons with size $N_u^H\times 1$, where $N_u^H$ is the number of DOF for velocity on the coarse grid. In the end, the coarse level features are downscaled/decoded to fine grid velocity output $N_u^h\times 1$ by a densely connected layer. An illustration of the network architecture is presented in Figure \ref{fig:vel_lcn}.

Furthermore, we construct the following constraint loss function in the deep learning process
 \begin{equation}\label{eq:loss_u}
\displaystyle{ \min_{\theta} } \quad  \frac{1}{N} \sum_{i=1}^N (  \frac{||u_{\text{pred, i}} - u_{\text{true},i} ||_2 }{ ||u_{\text{true},i} ||_2} + \ytrev{\beta}  || B (u_{\text{pred},i} - u_{\text{true}, i})||_2)
 \end{equation} 
where $ u_\text{pred} = \mathcal{N}(\kappa; \theta)$ or $ u_\text{pred} = \mathcal{N}(f; \theta)$ depends on the training samples, $N$ is the number of samples, and $B$ is the $B_h$ matrix in \eqref{eq:vel_mat_single}, \ytrev{$\beta$} is a regularization constant.

We remark that the second term in the loss function \eqref{eq:loss_u} helps to retain the local mass conservation of the predicted solution. \ytrev{The value of $\beta$ strikes a balance between minimizing the relative $L^2$ error of the true and predicted values of the velocity and forcing the predicted velocity solutions satisfying the physical constraint. We note that if the relative $L^2$ error of the true and predicted values of the velocity is small enough, the physical constraint may be achieved automatically. 
When this part of the loss decrease to some extent, the physical based loss comes to play. In our work, we naively use the grid search to select $\beta$. We did several experiments on the training data when $\beta =\{10^{-4}, 10^{-3}, 10^{-2}, 10^{-1}, 10^{0}\} \times \frac{\text{std}(||u_\text{true}||)}{\text{std}(\overline{|M_{\text{true}}|})}$, where $\frac{\text{std}(||u_\text{true}||)}{\text{std}(\overline{|M_{\text{true}}|})}$ is the factor between the scales of physical based loss and the mean squared loss. Then we use $10^{-3} \times \frac{\text{std}(||u_\text{true}||)}{\text{std}(M_{\text{true}})}$ which gives the relatively better results in the experiments.  }

With this custom loss function, the training is then carried out using the Adam optimization algorithm to minimize the loss evaluated on batches of data from the training set.

\subsection{Transport problem using sparse learning} \label{sec:sat_single}
\subsubsection{Feed-forward map in transport problem}
\begin{figure}[!hbt]
\centering
\includegraphics[scale=0.3]{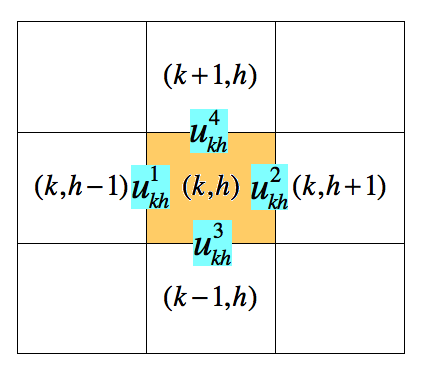}
\caption{ An illustration of $(k,h)$ cell and its neighborhood.}
\label{fig:5cells}
\end{figure}

Next, we will move to the saturation equation. The saturation equation \eqref{eq:sat_linear} can be discretized and written in the matrix equation as follows
\begin{equation}\label{eq:sat_lin_mat}
S^{n+1}  = S^{n} + \text{d} t ( F(u) S^{n} + R)
\end{equation}
Suppose the computational domain is partitioned using squares as fine cells. Consider the fine cell on $k$-th row and $h$-th column, 
according to the labels in Figure \ref{fig:5cells}, we can rewrite the saturation feed-forward map for this fine degree of freedom as follows
\ytrev{
\begin{equation}\label{eq:sat_relu}
\begin{aligned}
S_{k,h}^{n+1}  &= S_{k,h}^{n} + \frac{\text{d} t}{|e|} (\text{Relu} (u^1_{kh} \cdot \pmb{n} ) - \text{Relu} (u^2_{kh} \cdot \pmb{n}) + \text{Relu} (u^3_{kh}\cdot \pmb{n}) -\text{Relu} (u^4_{kh}) \cdot \pmb{n}) S_{k,h}^{n} \\
& + \frac{\text{d} t}{|e|} ( \text{Relu} (-u^1_{kh}\cdot \pmb{n}) S_{k-1,h}^{n}  + \text{Relu} (-u^2_{kh}\cdot \pmb{n}) S_{k+1,h}^{n}  \\
& +\text{Relu} (-u^3_{kh}\cdot \pmb{n})S_{k,h-1}^{n}  +\text{Relu} (-u^4_{kh}\cdot \pmb{n}) S_{k,h+1}^{n} )  + dt\, r_{kh}
\end{aligned}
\end{equation}
}
where $\pmb{n}$ is the unit outward normal vector, $|e|$ denotes the length of the edge.

\subsubsection{Network structure }
Given the flow velocity (either is trained from the deep network as described in the previous section, or obtained directly from MFEM solver), we can design a residual network to learn the dynamics of saturation equation. That is, we would like to approximate the feed-forward map from $S^n$ (input) to $S^{n+1}$ (output) using a deep neural network $\mathcal{M}$.

It is clear that the velocity field plays an important rule in the dynamics. However, if we directly let both velocity and saturation be input nodes to the network which are connected to the next layer neurons, the connections will be substantially large and may lead to an over-complicated neural network with redundancy. However, we can see from \eqref{eq:sat_lin_mat} that, the impact of velocity are stored in the matrix $F(u)$. Since the matrix $F(u)$ in \eqref{eq:sat_lin_mat} is sparse and the sparsity structure (a pentadiagonal matrix) is known, we will utilize the velocity values as well as the sparse structure of $F(u)$ to design some sparse weight matrices in the neural network. The sparse property will reduce the number of trainable parameters during the training.

To be specific, we first save the velocity solution obtained in Section \ref{sec:vel_single} in four directions separately (corresponding to four edges of each fine cell), i.e. $U = [u^1, u^2, u^3, u^4]$, where each column vector $u^{\cdot}$ has dimension $N_s^h \times 1$, where $N_s^h$ is the number of fine scale cells in the computational domain.

Next, we design a sparse custom layer, called \textit{Sparse Velocity Layer}. Motivated by Equation \eqref{eq:sat_relu}, we introduce four pentadiagonal matrices $W_i$ ($i=1,2,3,4$). Let $I$ and $J$ be the row and column indices of the pentadiagonal matrix, and $V_i$ be vectors consisting some random values drawn from a normal distribution, for $i=1,2,3,4$. The size of $I$, $J$, $V_i$ should be consistent. We initialize four sparse weight matrices as follows
\begin{equation}\label{eq:weight_def}
W_i = \text{sparse} (I, J, V_i), \quad i = 1,2,3,4.
\end{equation}
where $V_i$ are trainable parameters.

The Sparse Velocity Layer is defined by the following function
\begin{equation}\label{eq:self_sv_layer}
\phi_i(S^n, u^i) = \sigma \displaystyle( ( W_i \circ u^i ) S^n \displaystyle)
\end{equation}
where $\circ$ denotes the element-wise multiplication with broadcasting. To be specific, $W_i$ is a sparse matrix whose corresponding dense shape is $N_s^h \times N_s^h$, $u^i$ is an $N_s^h \times 1$ column vector, $W_i \circ u^i$ computes the element-wise product of each column of sparse matrix $W_i$ and $u^i$.  The latter multiplication of $( W_i \circ u^i )$ and $S^n$ is the standard matrix-vector multiplication.

\begin{figure}[!hbt]
\centering
\includegraphics[scale=0.35]{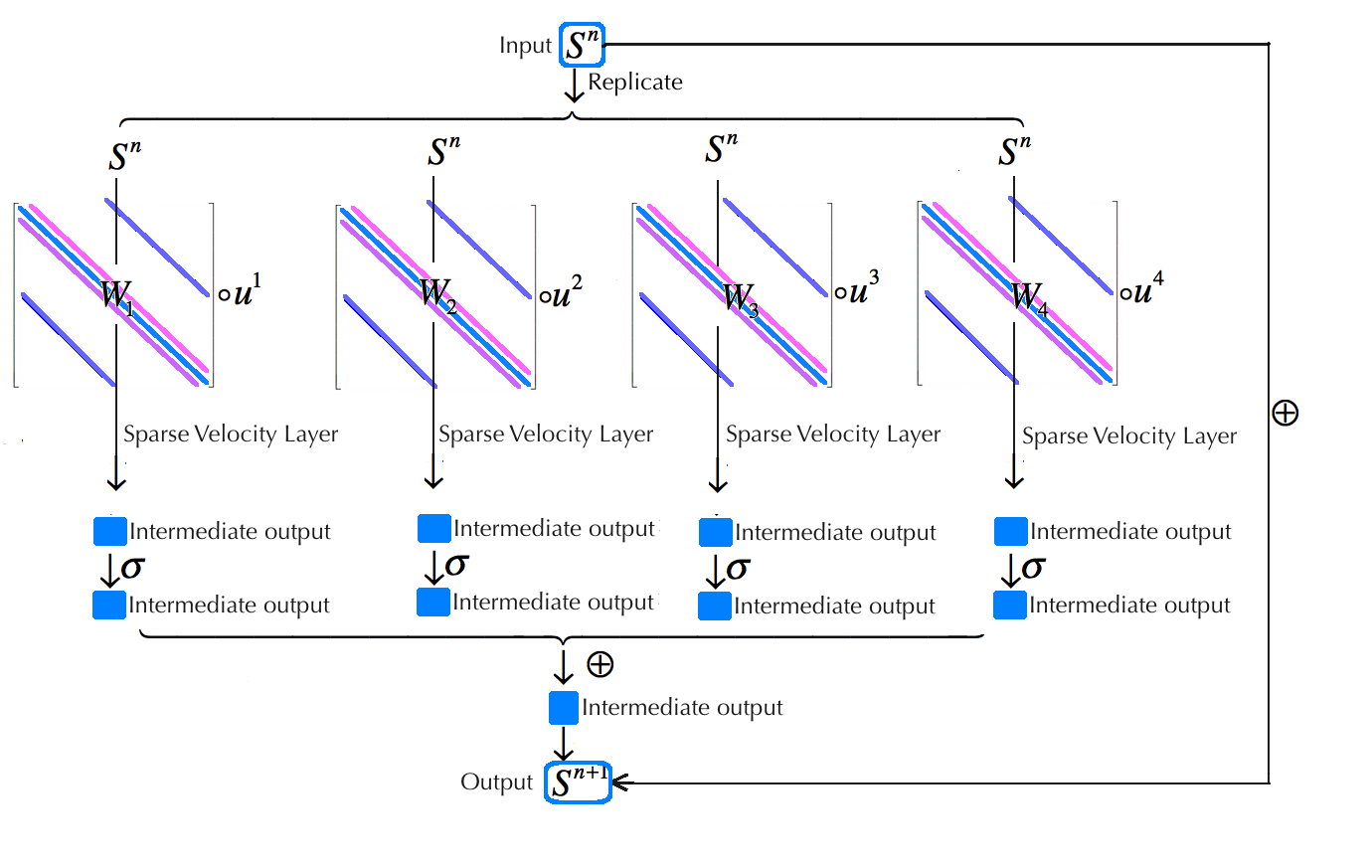}
\caption{An illustration of the neural network architecture for learning saturation problem.}
\label{fig:sat_resnet}
\end{figure}

Then we can design the following network to model the dynamics of the transport equation
\begin{equation}
\mathcal{M}(S^n; u) =   \sum_{i=1}^4 \phi_i (S^n, u^i)  + S^n
\end{equation}
where the network contains an addition of residual part $\sum_{i=1}^4 \phi_i (S^n, u^i)$ with the input $S^n$, as shown in Figure \ref{fig:sat_resnet}.

Once the network $\mathcal{M}$ is trained as an accurate surrogate model to approximate the map between saturation solutions at two consecutive time steps, we can use it iteratively several times to predict saturation solutions in the long range.

\ytrev{Remark 1: Remark 1: We remark that the equations of the problem are assumed to be known in this work, and we use both simulation data and the knowledge of the equations for the training of neural networks. Our approach doesn't require large data samples for the training, and can be used as an efficient surrogate solver for the investigated problems. When there are only simulation data without knowing physics, one may train surrogate models as a supervised learning task using methods like Gaussion process for small dataset, and some deep learning approaches for large data set. These data-driven methods usually relies on sufficient training data, which may be expensive to obtain for complex problems. If the underlying physics are also known, there are extensive work on physics informed neural networks (PINN) \cite{PINN1,PINN2}. This approach naturally encodes any underlying physical laws as prior information, and is successfully applied to solve PDE problems constrained to obey the law of physics that govern the data. Another development of physics-constrained surrogate modeling for stochastic PDEs without simulation data are proposed in \cite{PCDL_nz}. This work incorporates the governing equations in the loss function to enforce physical constraints without solving PDEs to obtain labeled/output data for the training. These developments shed lights on the combination of physics constraint surrogates and data driven surrogates for PDE systems. 
Our future work includes investigating the problems when there are uncertainty in the underlying physics or only part of the physics are known, and developing methods to combine physics and data based on current approach.}

 \subsection{Numerical results for the single phase flow case}
 \subsubsection{Numerical test on neural network approximation of the flow equation } 

\begin{figure}[!h]
\centering
\includegraphics[scale=0.27]{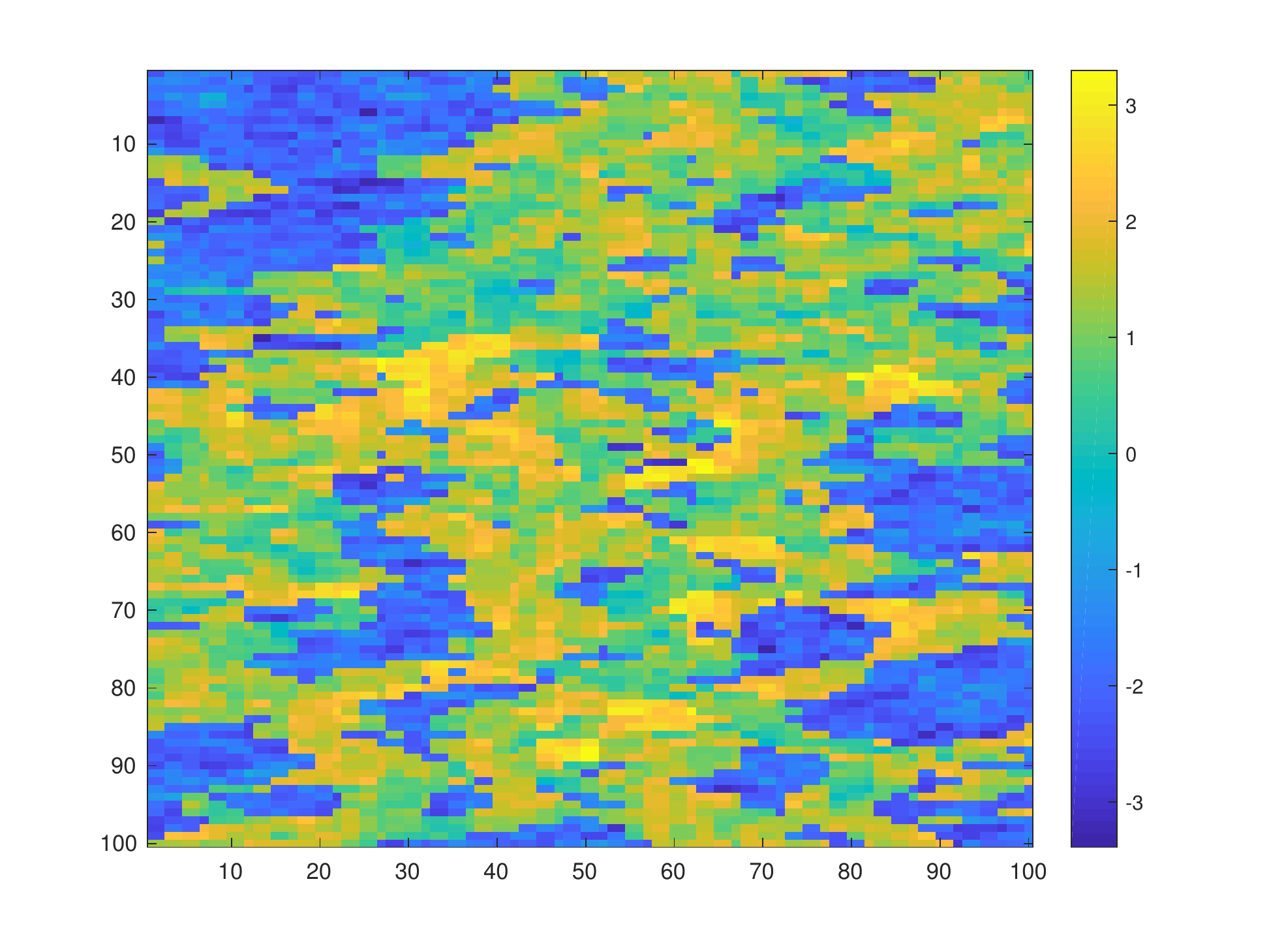}
\includegraphics[scale=0.27]{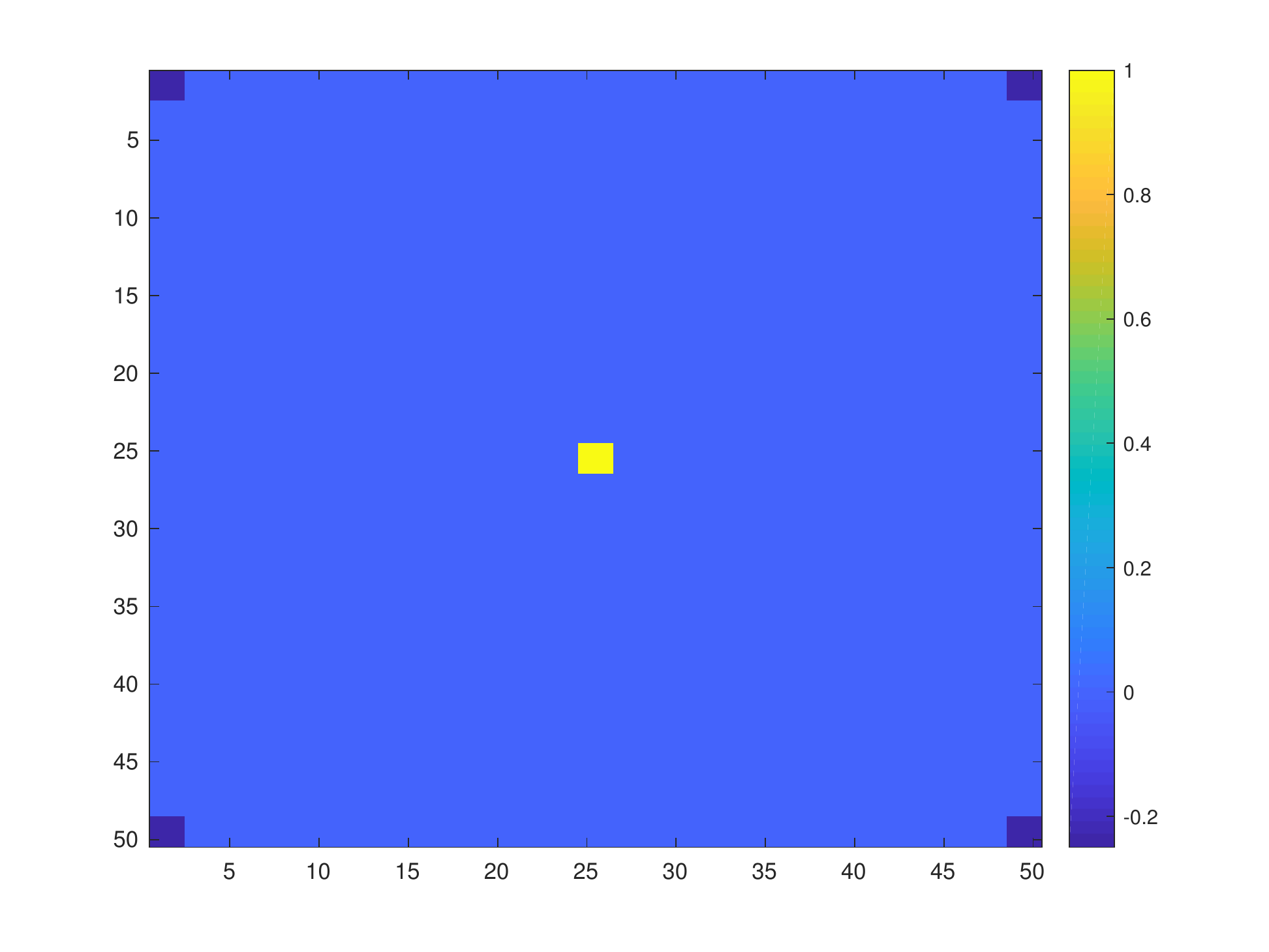}
\includegraphics[scale=0.27]{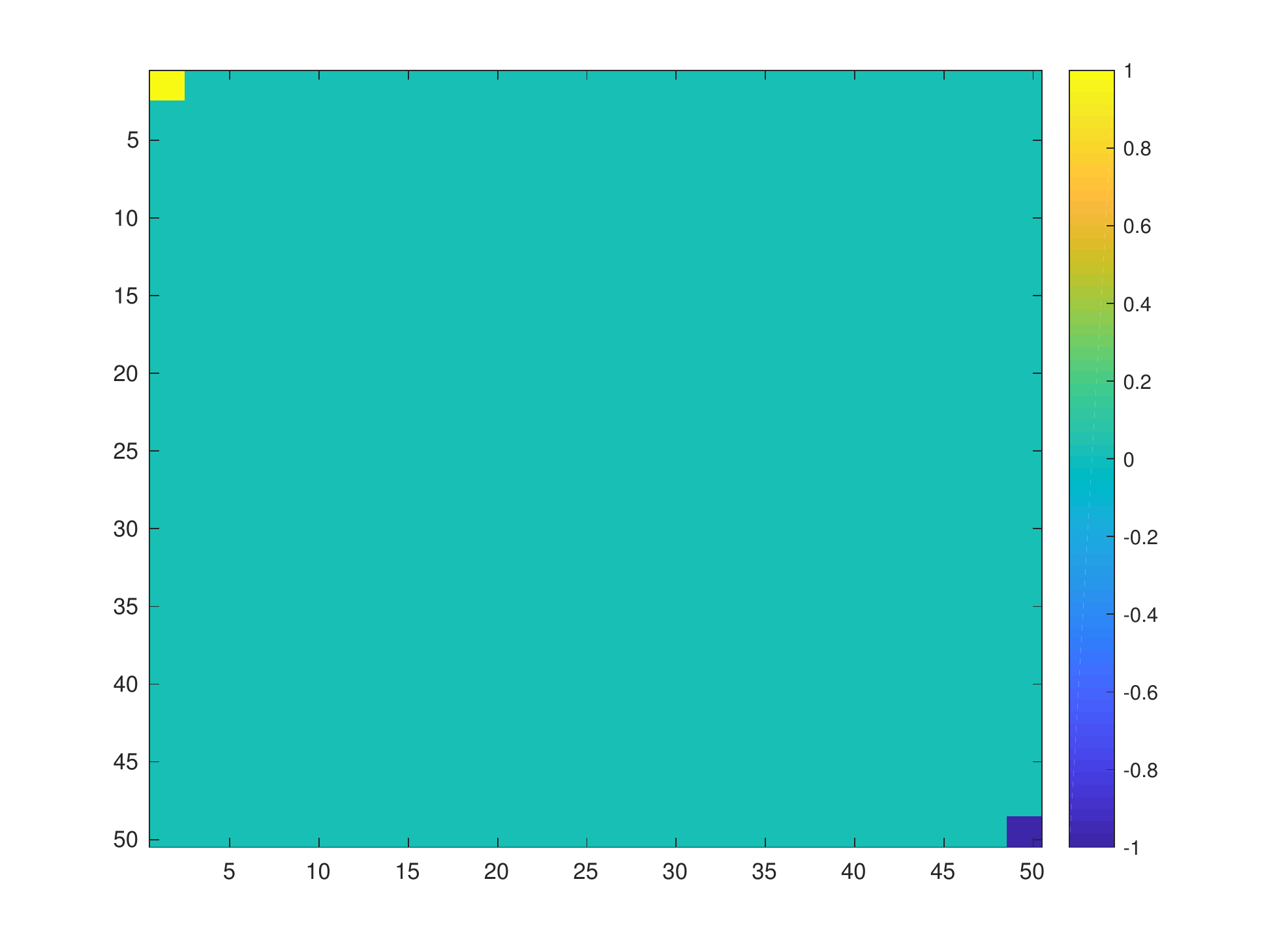}
\caption{ Left: Log scale of the permeability field (SPE10 model layer 53). Middle: Illustration of five-spot source in single phase flow problems. Right: Two-well source in two phase flow problems.}
\label{fig:perm}
\end{figure}

In this example, we will validate our proposed network for single phase flow problem. We consider the five spot reservoir configuration of the source term. Four injection wells are placed at the corners of the computational domain, and one production well is located at the center of the domain. By randomly varying the values of the injection rates, we generate different sources terms. Using these samples, we apply the mixed fem to solve for the velocity solutions. The permeability field (shown in Figure \ref{fig:perm}) remains the same among samples. The corresponding source and velocity pairs form the samples in the training process. In this example, we have $1000$ samples for training and $250$ for validation. The learning rate is chosen to be $0.008$. There are $100$ samples in each batch during the training. \ytrev{For training the network which approximate the flow problem, the number of epochs is $250$, and the total run time for the training is $29.2$ seconds when there are $1000$ samples.} After training, we take new source input and predict the velocity. We present, in Tables \ref{tab:vel_single_3network} and \ref{tab:vel_single}, the average errors between the predicted and true velocity solutions among $250$ test samples. We remark that all the network training in this work are performed using the Python deep learning API Keras \cite{chollet2015keras} with TensorFlow framework. 

We first compare the performance of (1) the proposed network as shown in Figure \ref{fig:vel_lcn}, denoted by locally connected network (LCN), where $N_p^H = 100$, $\alpha=4$, $\alpha'=8$, and the stride size is $4\times 4$; (2) a convolutional neural network, denoted by CNN, where we replace the locally connected layers with convolutional layers and all the other hyperparameters stay the same; (3) a fully connected network, denoted by DNN, where the network has the same number of layers as the LCN/CNN, and each hidden layer has $500$ neurons. The mean errors among the same set of testing samples are shown in Table \ref{tab:vel_single_3network}.
\begin{table}[!htb]
\centering
  \begin{tabular}{|c  |c  | c  | c|}
  \hline
 &$\norm{u_{\text{pred}} - u_{\text{true}} }_{L^2}$ (\%) &$\norm{u_{\text{pred}} - u_{\text{true}} }_{L^2_{\kappa}}$(\%)  & number of trainable parameters \\  \hline
LCN  &0.6 &0.7 &2,536,772 \\  \hline
CNN &0.8 &1.0 &2,525,708 \\  \hline
DNN &2.1 &2.2 &4,346,220\\  \hline
  \end{tabular}
\caption{Comparison of the true velocity solution (obtained from standard solver) and predicted velocity solution (obtained from trained neural networks) for three different types of networks. Mean errors are computed among 250 testing samples. Second column: velocity solution relative error in standard $L^2$ norm; Third column:  velocity solution relative error in weighted $L^2$ norm, \ytrev{where $||u||_{L^2_{\kappa}} = \int_{\Omega} \kappa^{-1} |u|^2$}. }\label{tab:vel_single_3network}
\end{table}
One can observe that, the proposed network shows great potential to approximate the mapping from input source terms to output velocity solutions. Our proposed network architecture works better than the CNN (\ytrev{where we only change the locally connected layers in the designed network structure with standard convolutional layers}), and also outperforms the DNN network. Compare the trainable parameters in the three network, LCN has slightly more parameters than that in the CNN, but has much smaller numbers than that in the DNN. 

\ytrev{Remark 2: We compare the LCN and CNN and restrict them to the same designed architecture in Figure 1, and the numerical results show that under same architecture, replacing convolutional layers by locally connected layers can give us slightly better results. The only reason for the comparison is that, we think the locally connected layer is more suitable simply in our case. In this work, the neural network we constructed is motivated by ideas from multiscale model reduction algorithms. From this perspective, we chose locally connected layers since it has a different set of filters applied at each different patch of the input. This is the analogy to the multiscale method, where local features of the underlying heterogeneity can be extracted using multiple different basis function in different local coarse regions. The locally connected layers is a good choice and can provide sufficient learning capacity to extract rich hidden features. However, there are many well established neural networks for solving problems of the similar interest. CNNs with deeper layers were also investigated in \cite{CNN} to automatically extract multiscale features from high dimensional input. Moreover, the deep convolutional encoder-decoder networks designed \cite{deepconv_Zabaras} shows very good results as surrogate models to solve stochastic PDEs. The method, on the other hand, doesn't require any explicit intermediate dimension reduction method and provide uncertainty estimates. These work provide alternative approaches and gain significant progress in surrogate modeling and uncertainty quantification.}

Let the standard loss refer to the mean of relative weighted $L^2$ errors
\begin{equation}\label{eq:loss_u_mse}
{  \frac{1}{N} \sum_{i=1}^N  \frac{|| u_{\text{pred, i}} - u_{\text{true},i} ||_2 }{ ||u_{\text{true},i} ||_2} }
\end{equation}
where there is no additional constraint as in \eqref{eq:loss_u}.

We next investigate LCN with the standard loss function (see \eqref{eq:loss_u_mse}) and the proposed loss function (see \eqref{eq:loss_u}). We can observe from Table \ref{tab:vel_single} that, using the proposed loss function with physical constraint, the average relative $L^2$ errors are smaller than the ones using standard loss. For the last column in the table, $\overline{M^i_{\text{true}} - M^i_{\text{pred}}}$ measures the errors for the local mass, where we first compute the mean of the local mass difference between predicted solutions and true solutions among all fine cells, then the average among all testing samples. With the additional constraint in the loss function, the locally mass conservative property is enhanced. 

\ytrev{Remark 3: There are many efforts to enforce physical constraints in many deep learning tasks recently. For example, the deep fluids CNN architecture is proposed in \cite{deep_fluids} which can synthesize divergence-free fluid velocities by introducing a novel stream function based loss function. There are also physical constrained models, such as \cite{PINN1, PINN2} employ the automatic differentiation techniques to differentiate neural networks regarding the input coordinates to obey the law of physics.}

\begin{table}[!htb]
\centering
  \begin{tabular}{|c  |c  | c | c |}
  \hline
 &$\norm{u_{\text{pred}} - u_{\text{true}} }_{L^2}$ (\%) &$\norm{u_{\text{pred}} - u_{\text{true}} }_{L^2_{\kappa}}$(\%) & $\overline{M^i_{\text{true}} - M^i_{\text{pred}} } $ \\  \hline
Standard loss  &0.4 &1.0 &1.64e-8	\\  \hline
 Loss with constraint &0.4 &1.0 &1.64e-9	\\  \hline
  \end{tabular}
\caption{Comparison of the true velocity solution (obtained from standard solver) and predicted velocity solution (obtained from trained neural networks) using standard loss and loss with constraints. Mean errors among 250 testing samples. Second column: velocity solution relative error in standard $L^2$ norm; Third column:  velocity solution relative error in weighted $L^2$ norm; Fourth column: local mass relative error in $l^2$ norm. }\label{tab:vel_single}
\end{table}

We also present the comparisons when we take a different number of nodes $N_p^H$ as shown in the network \ref{fig:vel_lcn}. As shown in Table \ref{tab:vel_single_diff_neurons}, when we use more neurons in the intermediate layers, (which corresponds to the larger number of coarse grid properties $N_p^H$, for example, $N_p^H = 25$ corresponds to $5\times 5$ coarse features, and so on), the trained network produce more accurate predictions. However, as $N_p^H$ become larger, the prediction errors decrease slower. This is due to the fact that the more complex of the network, the harder to train or the more training samples it requires. In our experiments, we use the same number of training sample in all three cases, so the results are acceptable.
\begin{table}[!htb]
\centering
  \begin{tabular}{|c  |c  | c | c |}
  \hline
$N_p^H$ &$\norm{u_{\text{pred}} - u_{\text{true}} }_{L^2}$ (\%) &$\norm{u_{\text{pred}} - u_{\text{true}} }_{L^2_{\kappa}}$(\%) & $\overline(M^i_{\text{true}} - M^i_{\text{pred}} ) $ \\  \hline
25 	&16.8 &26.0 &2.1e-9 \\  \hline
100  &0.4 &1.0 &1.64e-9	\\  \hline
225 &0.6  &0.8 &1.0e-9 	\\  \hline
  \end{tabular}
\caption{Comparison of the true velocity solution and predicted velocity solution using a different number ($N_p^H$) of neurons. Mean errors among 250 testing samples. All other hyperparameters in the network are the same. Second column: velocity solution relative error in standard $L^2$ norm; Third column: velocity solution relative error in weighted $L^2$ norm; Fourth column: local mass relative error in $l^2$ norm.}\label{tab:vel_single_diff_neurons}
\end{table}

\ytrev{
Remark 4: The last layer performs a downscaling step from coarse grid velocity to fine grid velocity solution. One can replace it with a sparse layer, where the weight matrix is sparse and takes into account the local effects. Designing the sparse layer requires the information about coarse-to-fine degrees of freedom map. If the dimensions of the fine/coarse scale space change, then one may need to redesign the last layer's sparse connection. We did both use dense layer and sparse layer in the last layer, and get similar results. For example, when $N_p^H = 100$ and sparse connections are used in the last layer, we got the mean error $0.7\%$ (compared to $0.4\%$ when dense connections are used in the last layer) in standard $L^2$ norm, and $0.3\%$ (compared to $1.0\%$ when dense connections are used in the last layer) in weighted $L^2$ norm.
}

In the end, we remark on the efficiency of the neural network. One fine solve of the problem \eqref{eq:vel_mat_single} using the Matlab direct solver takes $0.05$ seconds, and a prediction step using the trained network takes only $0.001$ second. We remark in this work, all experiments in Matlab are done on Intel(R) Xeon(R) CPU E5-1650, and experiments in tensorflow are done on GeForce GTX 1080 Ti.

\subsubsection{Numerical test on neural network approximation of the saturation equation }
In this section, we show a numerical example using the proposed neural network to learn the saturation dynamics. We still use a specific five spot configuration for the source term, and first obtain the velocity solution as described in Section \ref{sec:vel_lnn}.
Next, the transport equation for saturation is solved using the finite volume scheme \eqref{eq:sat_lin_mat} given initial condition $S^0$. We choose the time step \ytrev{size} $\triangle t = 1$, and solve the problem $1200$ time steps to obtain a series of saturation solution on the fine grids, i.e. $[S^1, S^2, \cdots, S^{1200}]$. 

In order to train the dynamics between $S^n$ and $S^{n+1}$, we choose the solution pairs $(S^i, S^{i+1})$, $i=1,\cdots, N_{\text{train}}$ as training samples, and $(S^i, S^{i+1})$, \ytrev{$i=N_{\text{train}}+1,\cdots, 1200$} for validation. Then we use \ytrev{$S^{N_{\text{train}}+1}$} as a test sample for prediction. As mentioned in the previous section, once the feedforward map is trained, we can use it multiple times to predict solutions at later times. Our goal in this example, is to apply the trained network $\mathcal{M}$ for $1200-N_{\text{train}}$ times, and compare the prediction with the true solution $S^{1200}$. \ytrev{During the training, the number of epochs is $500$, and the total run time for the training is $54$ seconds when there are $1000$ samples.}

We will first take $N_{\text{train}}=1000$, and test the performance of (a) the proposed network with the introduced Sparse Velocity Layer as in \eqref{eq:self_sv_layer}, (b) Sparse Connected Layer, where the weight matrices to be trained only have the same sparsity structure as in \eqref{eq:weight_def}, but are not multiplied by velocity, (c) densely connected layers. The results are presented in Table \ref{tab:sat_lin}. We remark that, the other hyperparameters, such as batch size, learning rate, training epochs, etc, in (a)-(c) are chosen as the same. The errors between the predicted solution $S_{\text{pred}}$ and true solution $S_{\text{true}}$ are measured by
\[
E_s  = \frac{\norm{S_{\text{pred}} - S_{\text{true}} }_{L^2}  }{\norm{S_{\text{true}} }_{L^2}}.
\]
\begin{table}[!h]
	\begin{center}
		\begin{tabular}{|c  |c  | c | c |}
			\hline
			&Sparse velocity layer &Sparsely connected layer &Densely connected layer \\  \hline
			$E_s(1)$  (\%)  &0.008   &0.02    &2.67	\\  \hline
			$E_s(10)$  (\%)  &0.08  &0.21 &9.56 \\  \hline
			$E_s(100)$  (\%)  &0.91  &1.85  &blow up	\\  \hline
			$E_s(199)$  (\%)  &2.20 &4.16   &blow up		\\  \hline
			 Trainable paras \#  &59,200 &59,200 &6,252,500  \\ \hline
		\end{tabular}
	\end{center}
	\caption{$E_s(n)$ denotes the relative error after $n$ time steps. Errors between true (obtained from standard solver) and predicted saturation (obtained from trained neural network), single phase flow. Second column: using proposed Sparse Velocity Layer; Third Column:  using Sparsely Connected Layer; Fourth Column: using Densely Connected Layer.}\label{tab:sat_lin}
\end{table}

\begin{figure}[!h]
\centering
\includegraphics[scale=0.27]{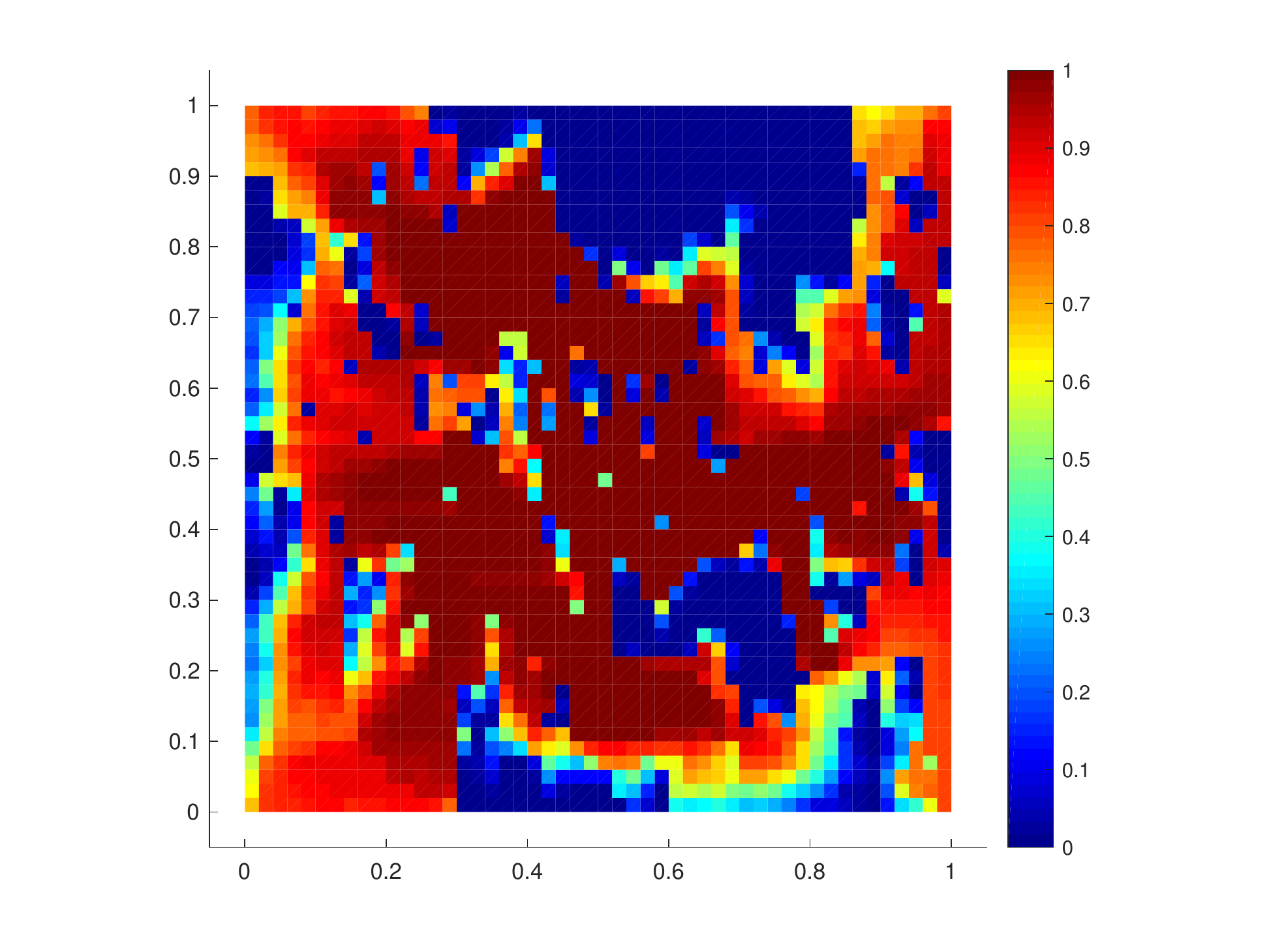}
\includegraphics[scale=0.27]{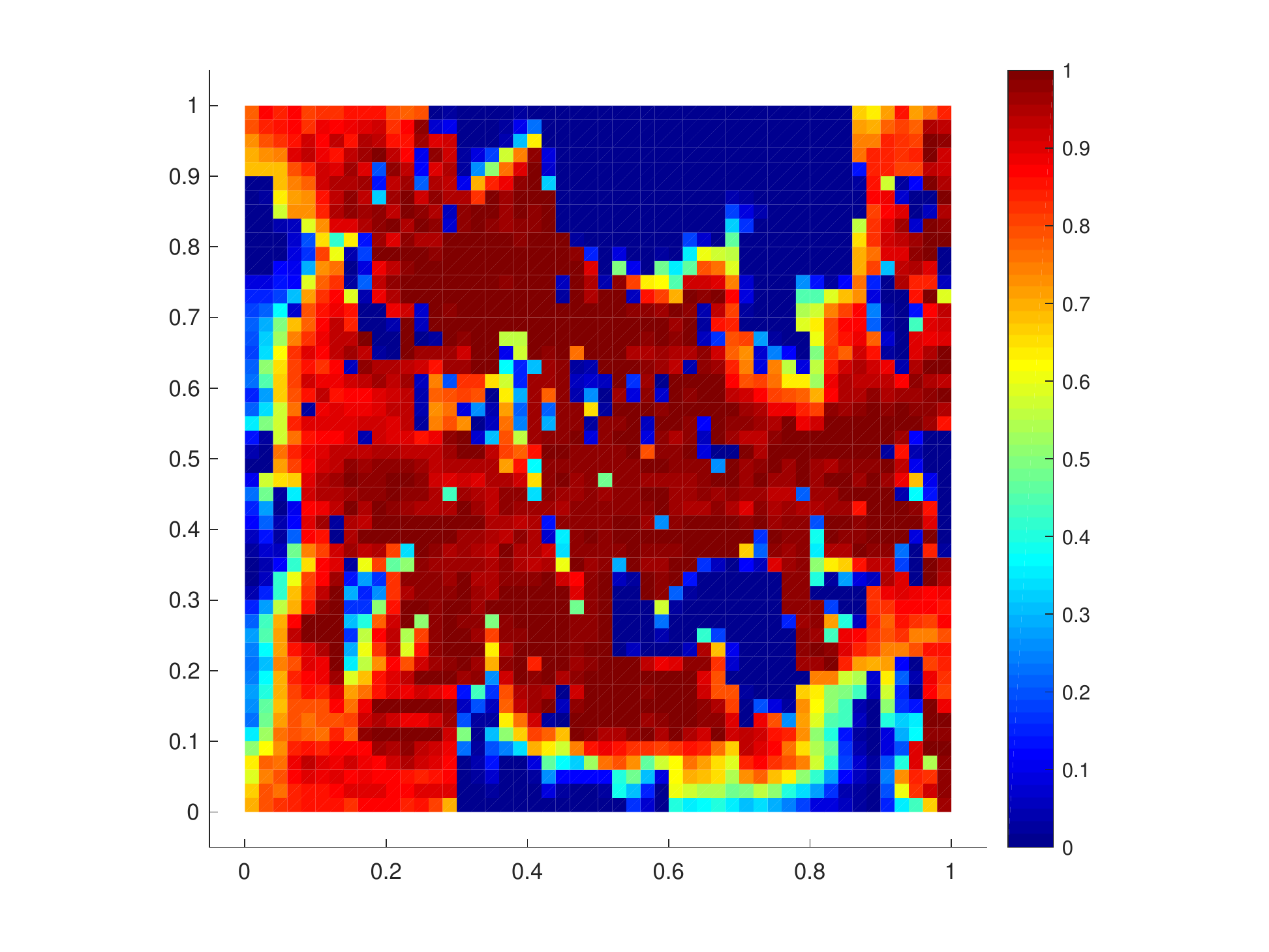}
\includegraphics[scale=0.27]{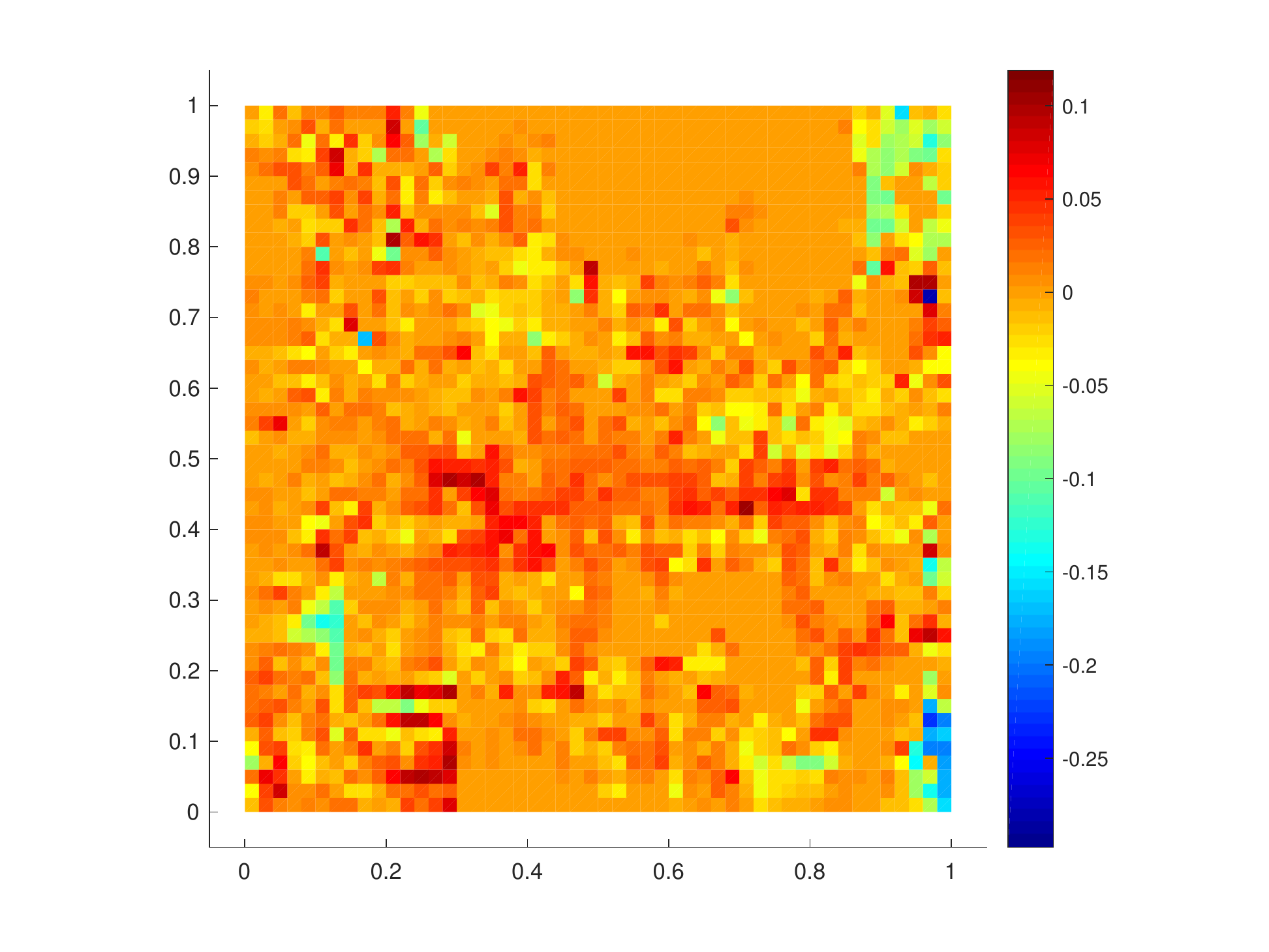}
\caption{Given solution at time step $1001$, use trained neural network to predict iteratively $199$ times. Comparison of saturation: left: true solution at time step $1200$, middle: predicted solution at time step $1200$, right: difference between the true and predicted solutions.}
\label{fig:sat_linear}
\end{figure}

From Table \ref{tab:sat_lin}, one can see that the number of trainable parameters is relatively small compared with a densely connected network using our proposed network. Moreover, the relative error of the predicted solution (the second column in the table) is very small even after $199$ time steps. However, using just sparsely connected layers without velocity information, the network prediction becomes a bit worse after some time (the third column in the table). With densely connected layers, the network works produce unreliable results (the fourth column in the table), this is due to a large number of the trainable parameters are hard to be trained effectively considering the limited number of training samples. What's more, this densely connected network didn't take into account the velocity field either. A comparison of the predicted solution using the proposed method and true solution in the long range is presented in Figure \ref{fig:sat_linear}, where we note the accuracy of the predicted solution.

Furthermore, we also test the performance of the proposed network by taking a different number of samples in the training set. The results are shown in Table \ref{tab:diff_samples_1ph}. It can be observed that our proposed method is quite robust. With more training samples, the network can be trained better. As we increase the predicted time steps, the errors grow reasonably. 


\begin{table}[!h]
\begin{center}
  \begin{tabular}{ |c  |c  | c | c |}
  \hline   
   \backslashbox{samples \#}{\# of time steps}    & 199 & 399 & 599	\\  \hline
1000  &2.20 & &	\\  \hline
 800 &2.65  &9.64 & 	\\  \hline
 600 & 5.11 &17.7 &68.4   \\  \hline
  \end{tabular}
  \end{center}
\caption{With different numbers of training samples and different predicted time steps, errors between predicted saturation (obtained from trained neural network) and true saturation (obtained from standard solver).} \label{tab:diff_samples_1ph}
\end{table}

In the end, we just want to mention that, in this single phase setting, using our proposed neural network $199$ times to predict the final solution takes only $0.466$ s. Solving the saturation equation using a direct solver in Matlab $199$ times takes $1.120$s.

\section{Two phase flow} \label{sec:two_phase}

Now, we consider incompressible two-phase flow. The gravity effects and capillary pressure will be neglected in the model. The flow equation reads
\begin{align}\label{eq:vel_2ph}
 u   = - \lambda(S) \kappa \nabla p \quad \quad &\text{in}  \quad D\\
\text{div}  (u) = r \quad \quad &\text{in}  \quad D\\
u\cdot n = 0 \quad \quad &\text{on}  \quad \partial D
\end{align}
where $\kappa$ is the absolute permeability and,
\[
\lambda(S) = \displaystyle{\frac{\kappa_{rw} (S)}{\mu_w} +  \frac{\kappa_{ro} (S)}{\mu_o}}
\]
is the total mobility, which depends on the saturation states. $\kappa_{rw}$, $\kappa_{ro}$ are the relative permeability, $\mu_w$ and $\mu_o$ are the viscosity, for water (denoted by subscript $w$) and oil (denoted by subscript $o$), respectively. Notice that, $\kappa_{rw}$, $\kappa_{ro}$ are nonlinear functions of $S$.

The saturation equation of $S_w$ for the water phase is given by
\begin{equation*}
\displaystyle{ \frac{\partial S_w}{\partial t} + u \cdot \nabla f_w(S_w) = r}
\end{equation*}
where $f(S_w) =\displaystyle{ \frac{\kappa_{rw}(S_w)/ \mu_w }{\kappa_{rw}(S_w)/\mu_w+\kappa_{ro}(S_w)/ \mu_o} }$. We simplify the notation and use $S$ for the saturation $S_w$ from now on.

The transport equation can also be solved using the finite volume method, and the time derivative is discretized using backward Euler. On a fine grid $K_i$, the value $S_i$ at time $t^{n+1}$ can be obtained by
\begin{equation} \label{eq:sat_2ph}
\displaystyle{ S_i^{n+1} = S_i^{n} + \frac{\text{d} t}{|K_i|}  [ -\sum_{e_j \in \partial K_i} F_{ij} (S^{n+1}) + f_w(S^{n+1})r_i^-  +  r_i^+ ] }
\end{equation}
\ytrev{where $r_i^- = \min(0,r_i)$, $r_i^+= \max(0,r_i)$. And} $F_{ij}$ is the upstream flux, i.e.
\begin{equation} \label{eq:upwind_flux_2ph}
F_{ij}(S^{n+1}) = \left\{
                \begin{array}{ll}
                  \int_{e_j} (u_{ij}^{n+1} \cdot n) f_w(S_i^{n+1}) \quad \text{ if }  \quad u_{ij}^{n+1} \cdot n  \geq  0\\
                    \int_{e_j} (u_{ij}^{n+1} \cdot n)  f_w(S_j^{n+1}) \quad\text{ if } \quad u_{ij}^{n+1} \cdot n  <  0
                \end{array}
              \right.
\end{equation}

Here $e_j$ denotes the edge between fine grid $K_i$ and $K_j$, $u_{ij}$ the velocity on the edge $e_j$.

For the two-phase problem, the flow and transport equations are solved sequentially follow the Algorithm \ref{alg:sequential}.

\begin{algorithm}
\caption{Sequentially solve flow and transport problem}\label{alg:sequential}
\begin{algorithmic}[1]
\Procedure{Two Phase}{$S^0, \text{d} t, T, t=0 $}\Comment{$S^0$: initial saturation, $\text{d} t$: time step, $T$: total time}
\State $S \gets S^0$
\While{$t < T$}
\State $\lambda(S^0) \gets \lambda(S)$
\State Solve $u$ from \eqref{eq:vel_2ph} using $\lambda(S^0)$
\State Solve $S$ from \eqref{eq:sat_2ph} using Newton-Ralphson method, $S^0$ and $u$
\State $S^0\gets  S$
\State $t \gets t+\text{d} t$
\EndWhile
\State \textbf{return} $S$
\EndProcedure
\end{algorithmic}
\end{algorithm}

One can observe from Algorithm \ref{alg:sequential} that, step $6$ requires Newton-Ralphson method which converges with a number of iterations. Also, a flow problem and a nonlinear transport problem are both solved at each time step. These will result in a heavy computational burden in practical problems. Our goal in this section is to construct powerful deep neural networks based on the experience in single-phase flow, to effectively approximate and enhance the computational efficiency of the coupled flow and transport solver.

\subsection{Flow velocity learning in two phase flow} \label{sec:vel_2ph}

In two-phase flow case, the source function $r$ is chosen to be a piecewise constant function. The source functions take nonzero values only at two fine cells, and zero values elsewhere. One of the two nonzero blocks is the fine cell at the upper left corner of the computational domain, the other is at the lower right corner of the computational domain. Considering compatibility, we assume one block has value $-Q$ acting as an injection well, the other has value $Q$ acting as a production well. The absolute permeability is highly heterogeneous and remains unchanged during the simulation, shown in Figure \ref{fig:perm}.

Given source terms and initial conditions, suppose we have solved the coupled flow and transport problem use Algorithm \ref{alg:sequential} for $n$ times with time step $\triangle t$. Let $\{r_1, \cdots, r_{n_s}\}$ be a set of source functions, \ytrev{where $n_s$ denotes the number of source samples}. For every source function $r_j$ ($j=1,\cdots, n_s$), we can obtain the velocity solutions $[u^1_j, \cdots u^{n}_j]$, and the saturation solutions are $[S^1_j, \cdots S^{n}_j]$. We will use the velocity and saturation solutions as samples to train neural networks. And we will omit the subscript $j$ for the simplification of notations. The design of neural networks in two-phase flow are based on the proposed methods in the single-phase flow, but some necessary modifications are made to handle the differences.

\subsubsection{Network Architecture}
For the two-phase system, we are interested in the map between the total mobility $\{\lambda(S), r\}$ and velocity field $u$ in the flow problem. As we can see from Algorithm \ref{alg:sequential}, in the coupled flow and transport system, the velocity fields are updated at each time step. The velocity updates come from the mobility change, and the mobility $\lambda(S)$ depends on the saturation solution from the previous time step. We remark that the absolute permeability $\kappa$ won't change during the simulation. However, the variation of $S$ changes the relative permeability.

Let \ytrev{$\big( \{\lambda(S^i), r^i\}, u^i \big)$} ($i = 1, \cdots n_t $) be the training samples, where $\{ \lambda(S^i), r \}$ is the input, and $u^i$ is the output, $n_t$ is the number of samples. We will follow the similar idea as in Figure \ref{fig:vel_lcn} to construct the neural network $\mathcal{N}_2$. However, there are some main differences we need to take care of. First, we have now both the mobility $\lambda(S)$ and $f$ as inputs. Compared with the single phase flow case, where we input $r$ as a single channel image, now we will place $\lambda(S)$ and $r$ in two channels of an image separately and use it as input. Second, we will add a few convolutional layers before the first average pooling layer. This is due to the complexity in $(\lambda(S^i), r^i)$ compared with the only source term $r$. The additional convolutional layers will extract important hidden features. After that, we will reduce the dimension by an average pooling layer. Then similar strategies are applied as described in the single-phase case \ref{fig:vel_lcn}.

\subsubsection{Numerical example for flow map in two-phase flow}
In the following, we present the results for two-phase flow velocity learning. Given a set of different source terms $\{ r_1, \cdots, r_{n_s}\}$, the coupled flow and transport system are solved $n=1000$ times with time step $\triangle t = 0.2$ for each source. Omitting the subscripts, we have $\big( (\lambda(S^i), r^i), u^i \big)$ ($i = 1, \cdots 1000\times n_s $) as training samples. After training, given a new source term $r_{\text{new}}$, we will take $\lambda(S_{\text{new}}^i), r_{\text{new}}^i)$ $(i=1, \cdots, 1000)$ are used for validation/prediction. 

In this example, the number of trainable parameters in this neural network is $2,625,236$. \ytrev{The number of epochs in the training is $1000$, and the total run time is $295$ seconds when there are $1000$ samples.} A prediction step of the trained neural network takes around $0.003$ seconds. Similar as in the single phase, local mass conservation is an important feature of the numerical velocity solution obtained from MFEM. This property ensures its adequacy in the transport simulations. We would like the solutions predicted by the proposed neural network also preserve this property. The constraint loss function is employed. By comparing the performance of the network when we use the constraint loss function \eqref{eq:loss_u} and \eqref{eq:loss_u_mse}, we notice that with constraint loss function, not only the local mass conservative property of the predicted solution is enhanced, the accuracy of the predicted solution in both $L^2$ and weighted $L^2$ norm is also improved. The results are shown in Table \ref{tab:vel_twoph}. 

\begin{table}[!htb]
\centering
  \begin{tabular}{|c  |c  | c | c |}
  \hline
Errors &$\norm{u_{\text{pred}} - u_{\text{true}} }_{L^2}$ (\%) &$\norm{u_{\text{pred}} - u_{\text{true}} }_{L^2_{\kappa}}$(\%) & $\overline(M^i_{\text{true}} - M^i_{\text{pred}} ) $ \\  \hline
 Loss with constraint &0.7 &1.3 &1.2e-8	\\  \hline
Standard loss  &0.8 &1.9 &2.2e-8	\\  \hline
  \end{tabular}
\caption{With different loss functions, mean errors  between the true velocity solution and predicted velocity solution in two phase flow. Second column: velocity solution relative error in standard $L^2$ norm; Third column:  velocity solution relative error in weighted $L^2$ norm; Fourth column: local mass relative error in $l^2$ norm.}\label{tab:vel_twoph}
\end{table}

\subsection{Saturation dynamics learning in two phase flow}
As we see from the saturation equation, \eqref{eq:sat_2ph} can be written in the following matrix equation
\ytrev{\begin{equation}\label{eq:sat_nonlin_mat}
S^{n+1}  = S^{n} + \text{d} t ( \tilde{F}(u^{n+1}, S^{n+1}) + R)
\end{equation}
}
where $\tilde{F}$ is the upwind flux, nolinearly depends on the saturation as described in \eqref{eq:upwind_flux_2ph}. Due to the nonlinearity, \eqref{eq:sat_nonlin_mat} needs to be solved using an iterative method. This makes it much more computationally expensive compared with the single flow problem.

To improve efficiency, a residual network $\mathcal{M}_2$ is expected to learn the dynamics of the saturation equation. Similar as before, the feed-forward map we are interested in is still from $S^n$ to $S^{n+1}$. The main difference in the two-phase flow lies in that, the velocity field driving the convection is time-dependent. During the learning, the velocity fields are no longer universal. We need to adjust the previously proposed network in Section \ref{sec:sat_single} to handle these issues.

\subsubsection{Network structure}

\begin{figure}[!h]
\centering
\includegraphics[scale=0.35]{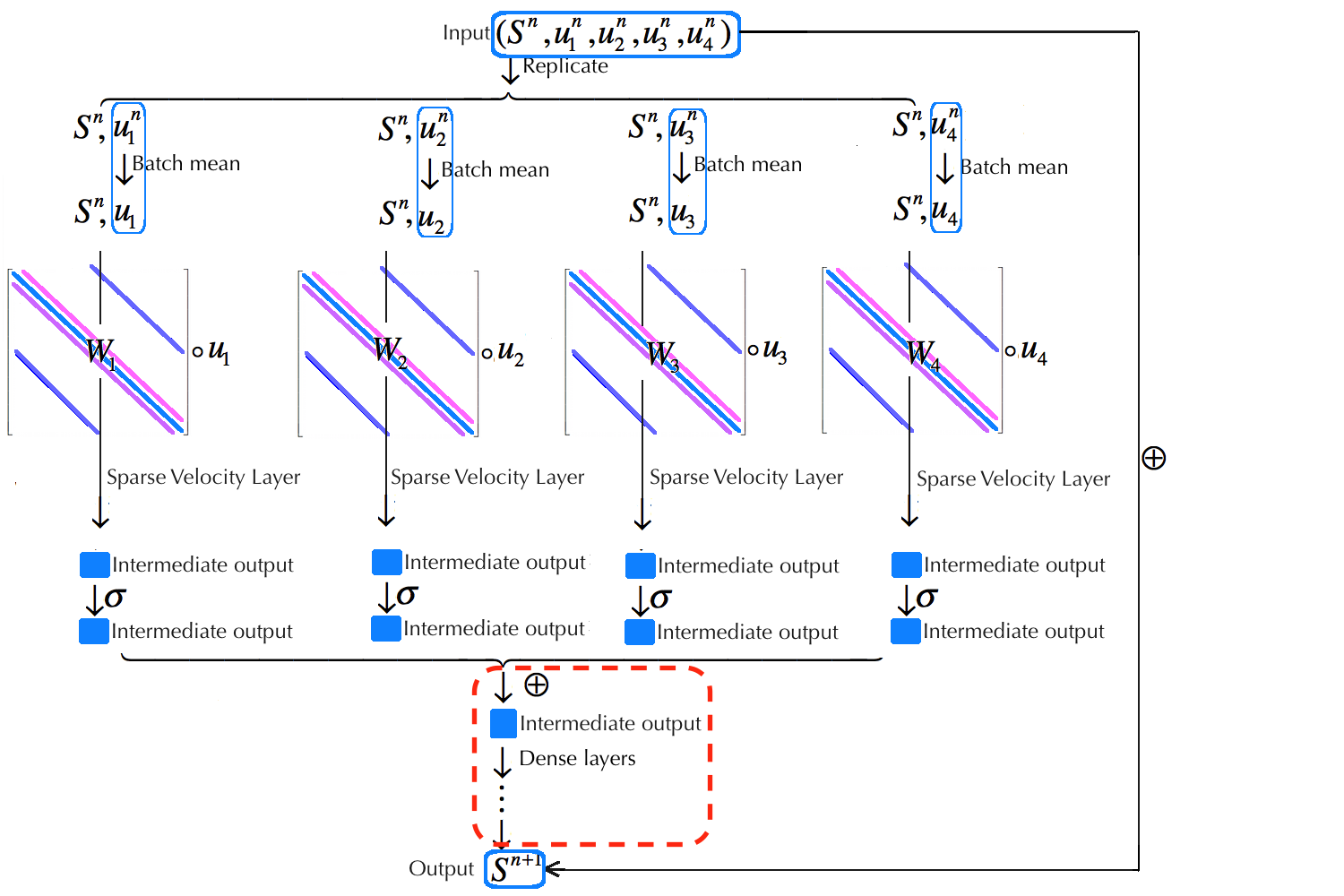}
\caption{An illustration of the neural network architecture for learning saturation problem in two phase flow.}
\label{fig:sat_resnet_2ph}
\end{figure}

We recall that, in the single-phase learning, we introduce the Sparse Velocity Layer to take into account the velocity impact on the transport problem, where a steady velocity field can be multiplied to the weight tensor element wisely in the layer construction. Different from the single-phase case, the velocity field is now changing at each time step. Thus, we need to take velocity fields corresponding to the saturation solution at each time step together as input to the network. However, it will introduce too many trainable parameters if we exactly use the velocity field corresponding to the saturation for each sample in the Sparse Velocity Layer. One workaround is to take the average velocity of each batch and use their mean in the Sparse Velocity Layer. Please see Figure \ref{fig:sat_resnet_2ph} as an illustration. Moreover, since the map we approximate is nonlinear, before the final residual add layer, we also add a few dense layers with nonlinear activation functions to capture the nonlinearity in the map we approximate.

\subsubsection{Numerical example of saturation learning in two-phase flow}
Here, we apply the proposed neural network \ytrev{shown in Figure} \ref{fig:sat_resnet_2ph} to learn the saturation dynamics in two-phase flow.

For a given set of source terms $\{r_1, \cdots, r_{n_s}\}$, we take the velocity and saturation solutions obtained from the coupled flow and transport solver, i.e. $[u_j^1, \cdots u_j^{1000}]$, and $[S_j^1, \cdots, S_j^{1000}]$ ($j=1, \cdots, n_s$).  We then separate the velocity in four directions, $u^i \rightarrow (u^i_1, u^i_2, u^i_3, u^i_4)$.  Abandon the subscript, we choose the states $(S^i, u^i_1, u^i_2, u^i_3, u^i_4, f^i)$ as training inputs,  and $S^{i+1}$ as training outputs, for $i=1,\cdots, 1000\times n_s$.  In the end, the velocity and saturation solutions at all time steps associated with a new source term $r_{\text{new}} \notin \{ r_1, \cdots, r_{n_s}\}$, are used for validation/prediction.

\begin{figure}[!hbt]
\centering
\includegraphics[scale=0.4]{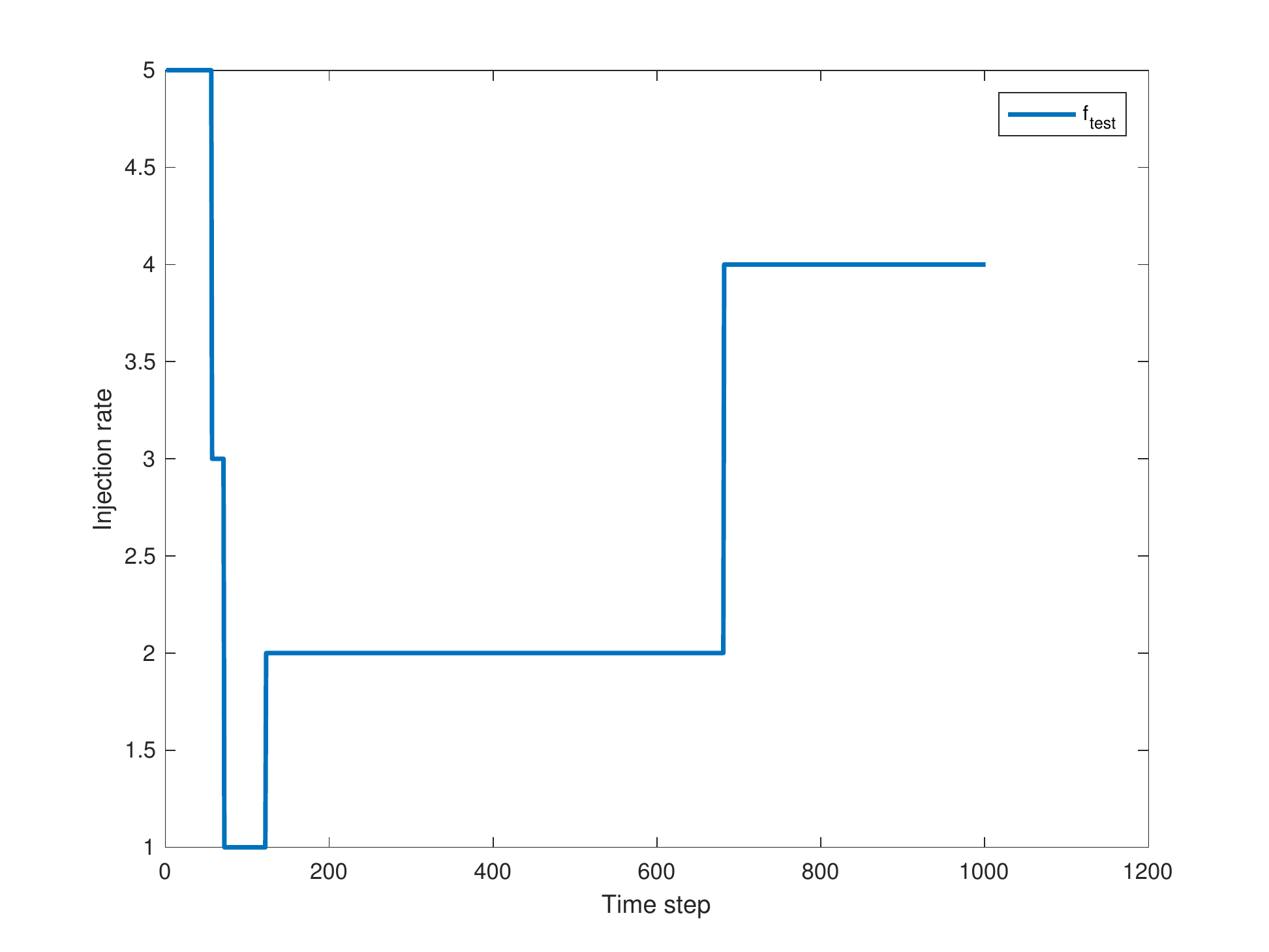}
\includegraphics[scale=0.24]{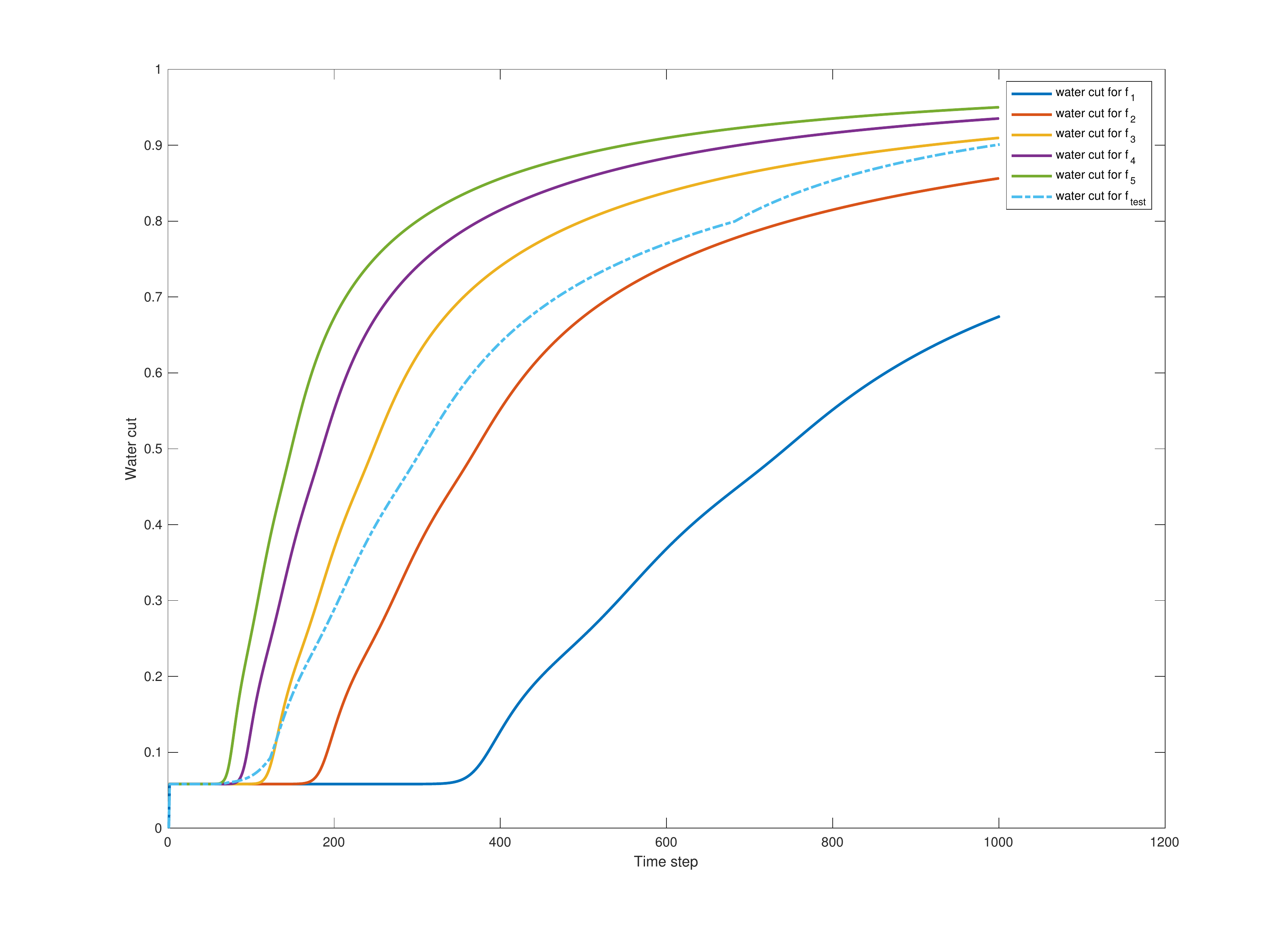}
\caption{ Left: a random source $r_{\text{test}}$ for generating test samples. Right: the water cuts at all time steps, for training sources $r_j$ ($j=1, \cdots, n_s=5$), and $r_{\text{test}}$. }
\label{fig:rand_ftest}
\end{figure}

In this example, for training, we choose $n_s = 5$, and $r_j=j$ at the injection well, for $j=1,\cdots n_s$. Here $r_j$ is a constant for all time steps. As for validation, we generate a random source as follows. First, the values of the source term is changing at four random time instants between $1$ and $1000$. Moreover, in each sub-interval, the values of the injection rate is also randomly chosen between $1$ and $5$. For example, we show a test source term in the left side of Figure \ref{fig:rand_ftest}. The corresponding water cuts for training source and testing source are shown in the right side of Figure \ref{fig:rand_ftest}. \ytrev{For training the network in the two phase flow, the number of epochs is $500$, and the total run time for the training is around $588$ seconds.}

We remark that the batch size in this proposed network may affect the accuracy of the trained network since we take the mean of velocity solution among each batch in the Sparse Velocity Layer.  We will test the performance of the proposed network (as shown in Figure \ref{fig:sat_resnet_2ph}) with different batch sizes.

\begin{table}[!htb]
\begin{center}
  \begin{tabular}{|c  | c |}
  \hline
 Batch size   &Mean error $E_s$ over 1000 test cases  (\%)\\  \hline
10 & 0.036    \\  \hline
100 & 0.039    \\  \hline
200 &0.064   \\  \hline
  \end{tabular}
  \end{center}
\caption{Two phase flow: For different batch size, mean errors between the true saturation solutions at time $n+1$ and predicted saturation solutions at time $n+1$ (predicted from saturation at time $n$). The test source term is shown in the left side of Figure \ref{fig:rand_ftest}. }\label{tab:sat_only}
\end{table}

Our numerical tests show that the errors are not increasing too much when the batch sizes become large. One reason for this may be the variations in the velocity field are not that dramatic. Although the mobility $\lambda(S)$ in each time step is different, the absolute permeability $\kappa$ still have a strong impact on the velocity solution. But the absolute permeability stays the same throughout the simulation in our experiment.

\subsection{Coupled flow and transport prediction using trained neural networks}

We have obtained two surrogate models using deep learning, namely, $\mathcal{N}_2$ for the velocity, and $\mathcal{M}_2$ for the saturation equation in two-phase flow case. They have shown their efficiency and accuracy separately in the previous sections.

Similar as the idea in the single phase flow case, once the feed-forward map for saturation equation is trained, we can use it multiple times to predict solutions at later times. However, in the two phase flow case, we also need to feed in velocity field at each prediction step. That is not a problem for us, since the velocity field can be obtained using the network $\mathcal{N}_2$ we trained in Section \ref{sec:vel_2ph}. That is, we will use both of the trained network $\mathcal{N}_2$ and $\mathcal{M}_2$ iteratively to predict future solutions. The algorithm is as follows:

\begin{algorithm}
\caption{Deep learning for flow and transport problem}\label{alg:nn_2ph_alg}
\begin{algorithmic}[1]
\Procedure{Two Phase}{$S^n$, $\mathcal{N}_2$, $\mathcal{M}_2$, $m$}\Comment{ $S^n$: saturation at time step $n$, $N_u$ : network for velocity, $N_s$ : network for saturation, $m$: number of time steps for prediction}
\State $S \gets S^n$
\While{$i < m$}
\State $\lambda(S^n) \gets \lambda(S)$
\State Predict $u^n$ using $\mathcal{N}_2$ with input $\lambda(S^n)$
\State Predict $S^{n+1}$ using $\mathcal{M}_2$, with input $S^n$, $u^n$
\State $S \gets S^{n+1}$
\State $i \gets i+1$
\EndWhile
\State \textbf{return} $S$
\EndProcedure
\end{algorithmic}
\end{algorithm}

\subsubsection{Numerical example for applying trained networks to predict coupled system in the long range}

Our last numerical example shows the power of the trained networks in two-phase flow. For a new source term $r_{\text{test}}$ (shown in left side of Figure \ref{fig:rand_ftest})), we only need to take the solution $S^{1}$ as an input to the trained neural network. Then we will apply the networks multiple times to predict the saturation solution in the long range according to Algorithm \ref{alg:nn_2ph_alg}. 
Our goal in this example is to apply Algorithm \ref{alg:nn_2ph_alg} with trained networks $\mathcal{N}_2$ and $\mathcal{M}_2$ with $m=999$ as the number of predicted time steps, to predict the saturation $S^{1000}$. Then we will compare the predictions with true solutions obtained from Algorithm \ref{alg:sequential}. From Table \ref{tab:sat_coupled_2ph}, we observe that, the predictions are very close to the true solution. We have tested over many different cases of $r_{\text{test}}$ and got similar results, here we just take a few for illustration. 

\begin{figure}[!h]
\centering
\includegraphics[scale=0.27]{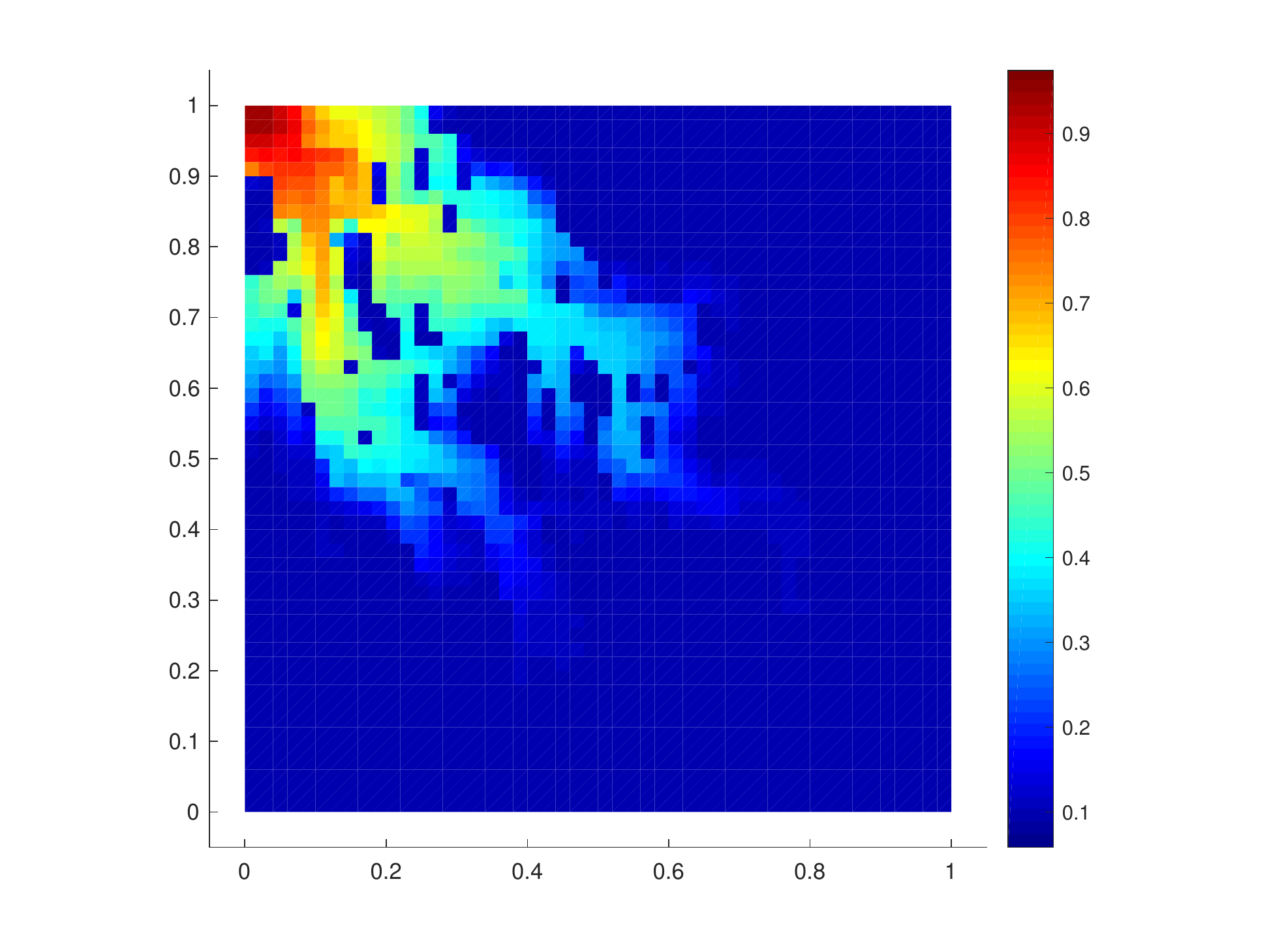}
\includegraphics[scale=0.27]{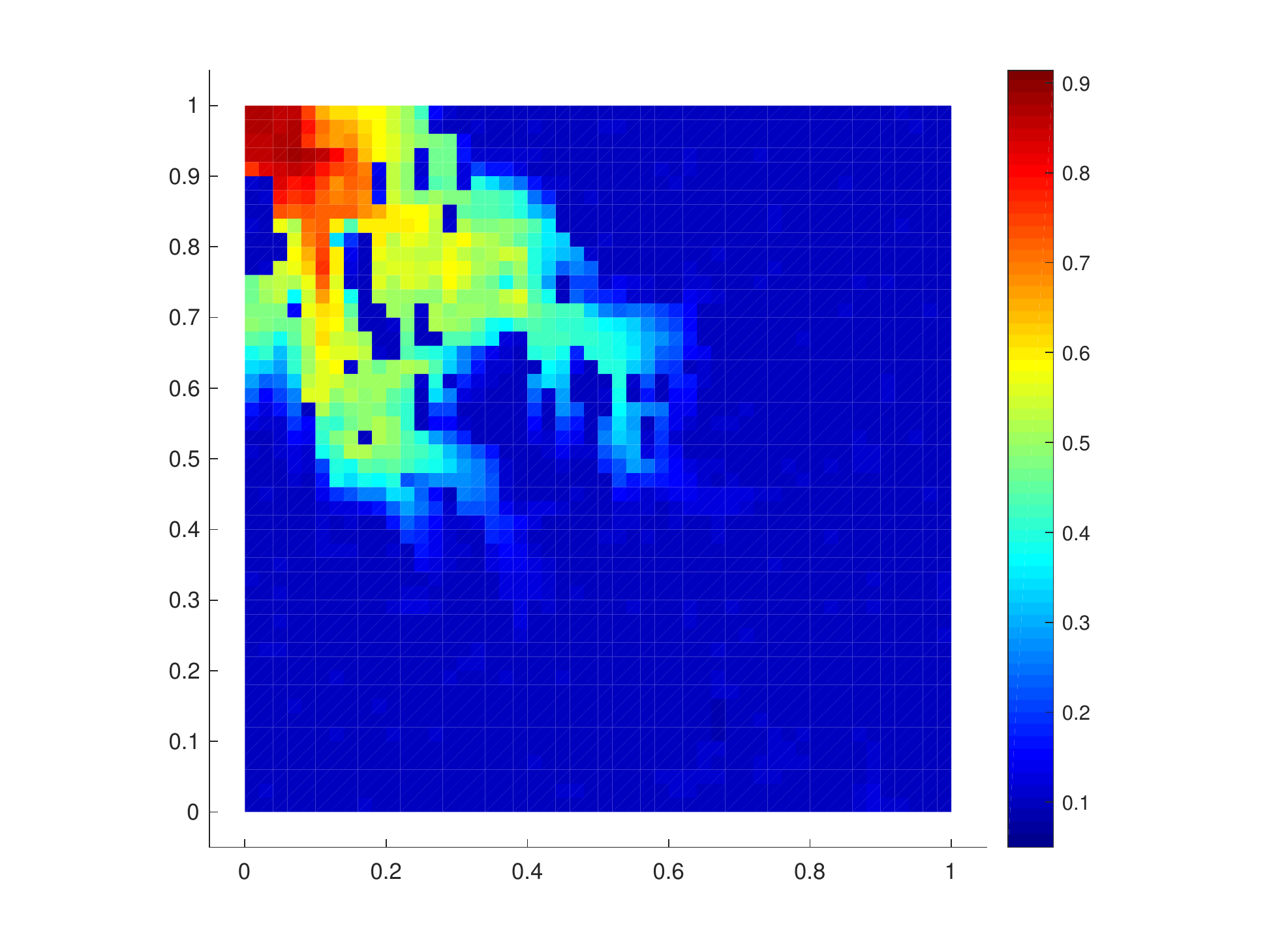}
\includegraphics[scale=0.27]{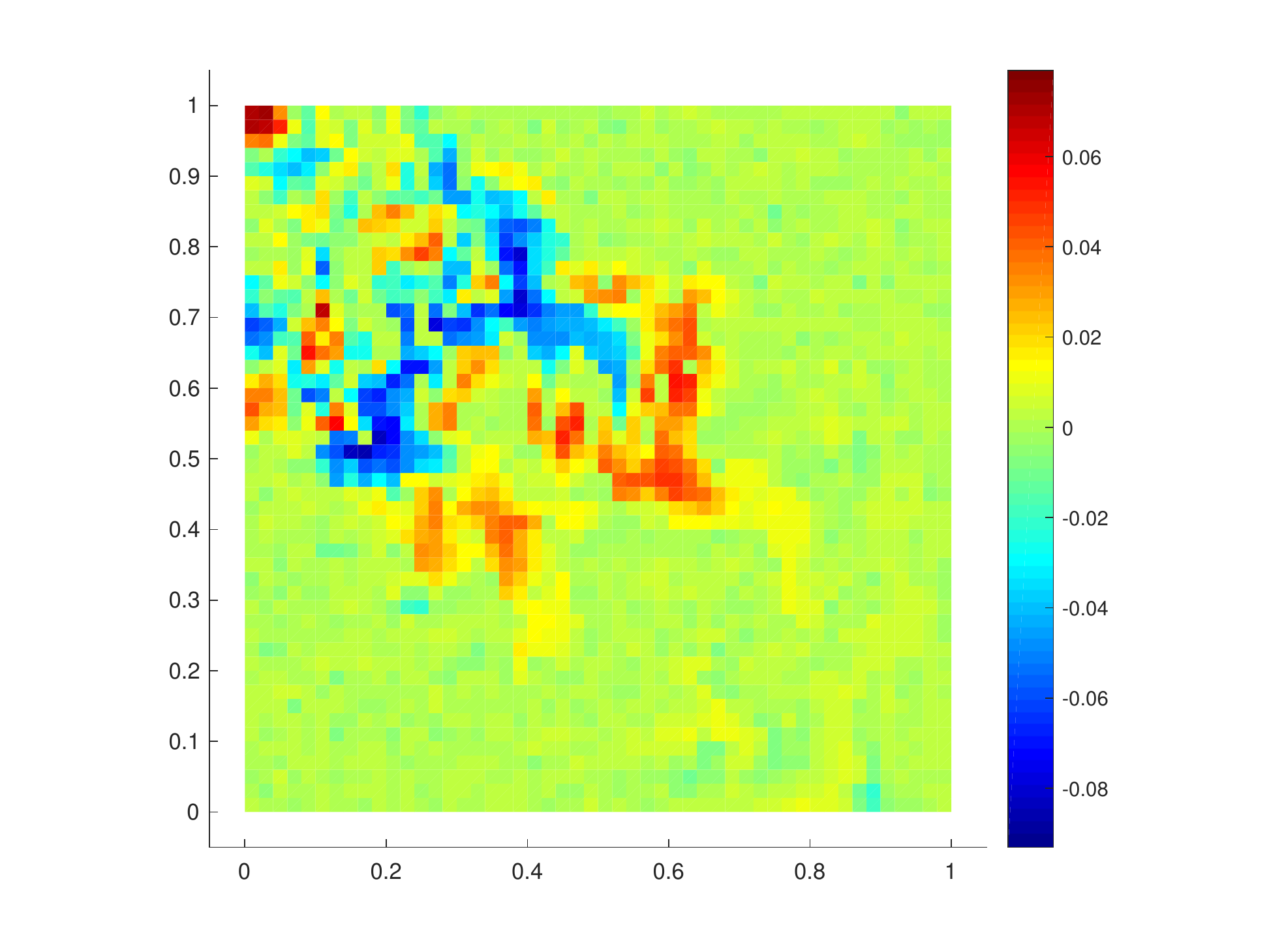}

\includegraphics[scale=0.27]{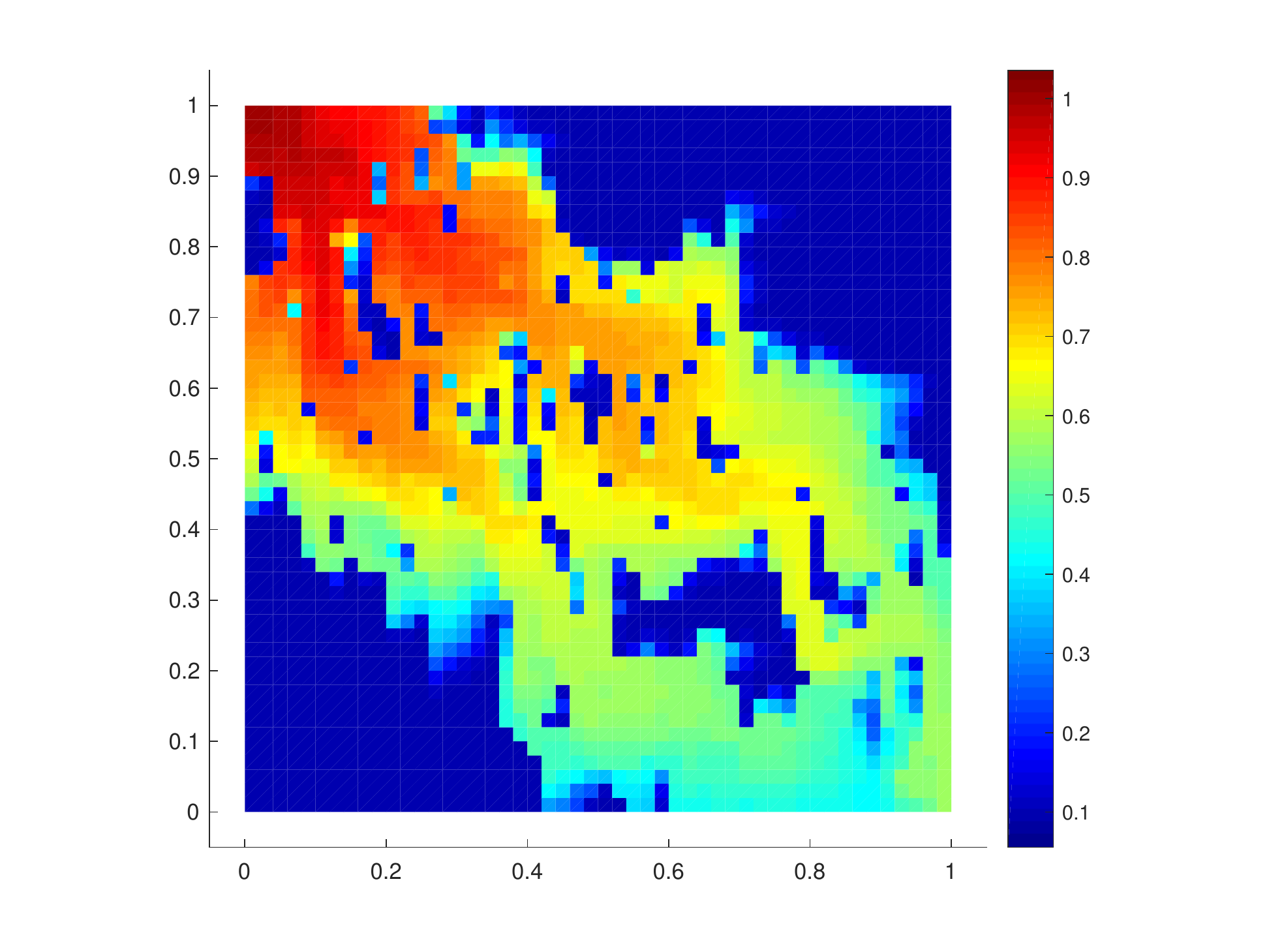}
\includegraphics[scale=0.27]{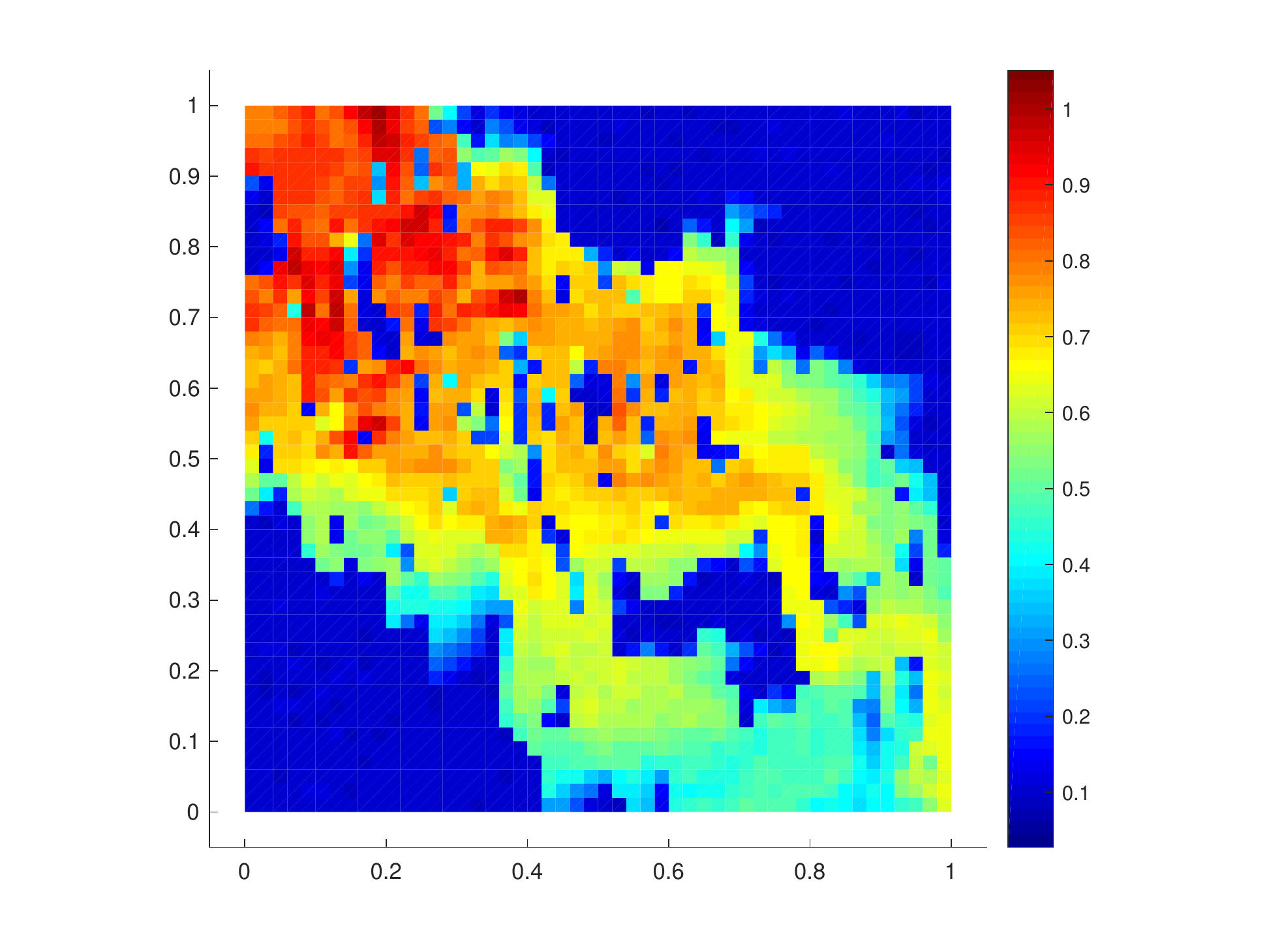}
\includegraphics[scale=0.27]{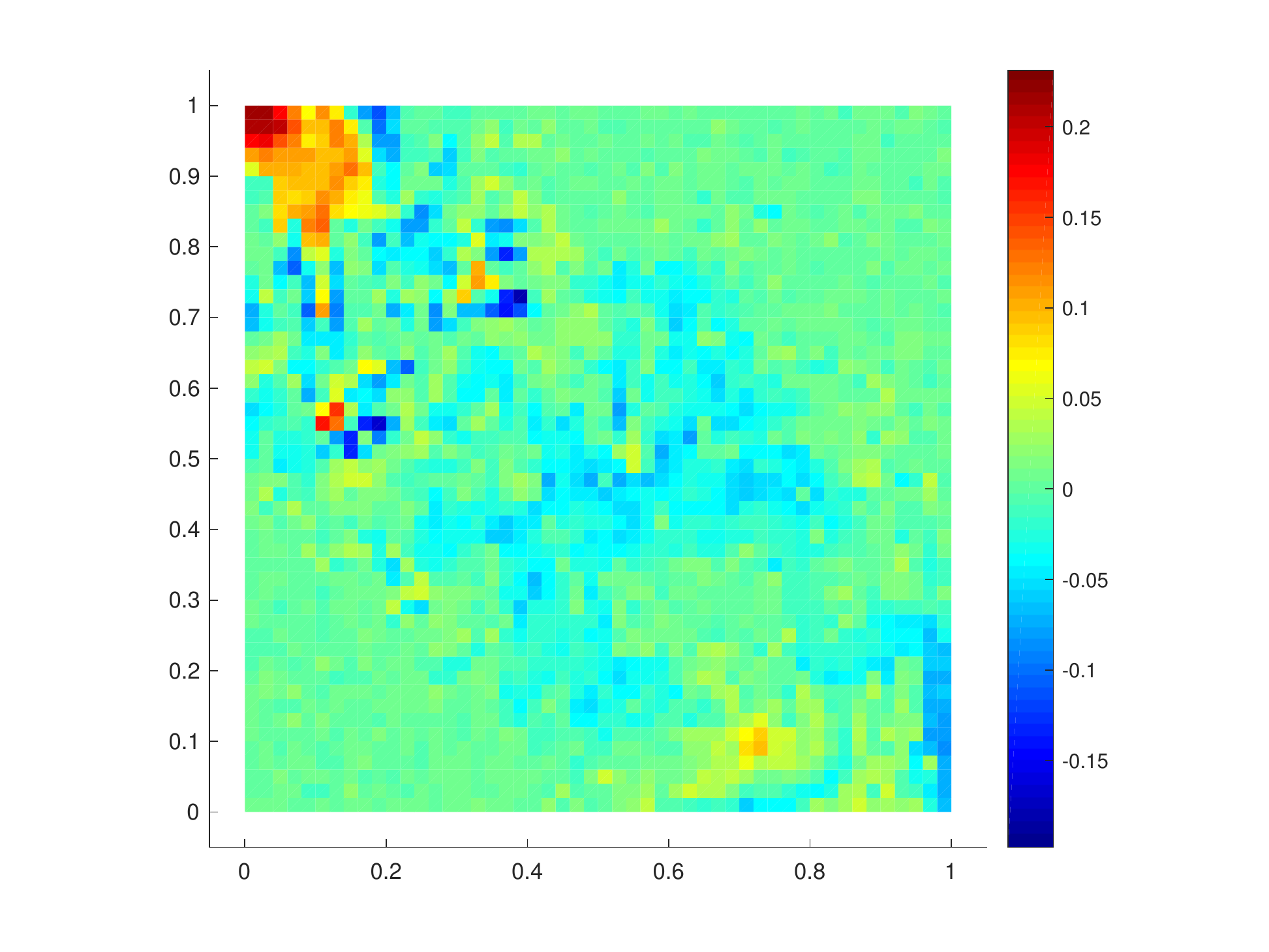}
\caption{ Comparison of saturation using Algorithm \ref{alg:nn_2ph_alg}. Given the initial solution, top row: after 50 prediction step, bottom row: after 950 predicted time steps. Left column: true solution, middle column: predicted solution, right column: difference between the true and predicted solutions. After 50 prediction step, the relative $L^2$ error between predicted and true solution is $7.0 \%$. After 199 predicted time steps, the relative $L^2$ error between predicted and true solution is $6.5 \%$}
\label{fig:sat_nonlinear}
\end{figure}

\begin{figure}[!hbt]
\centering
\includegraphics[scale=0.35]{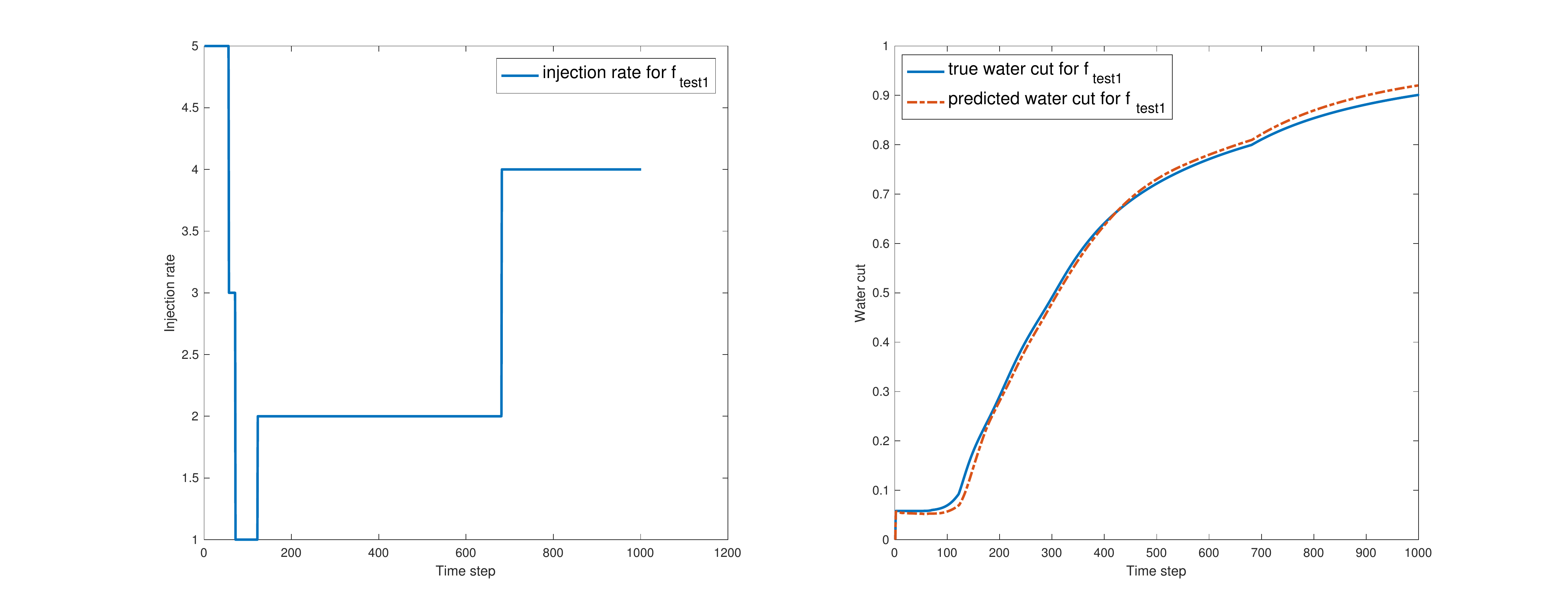}\\
\includegraphics[scale=0.35]{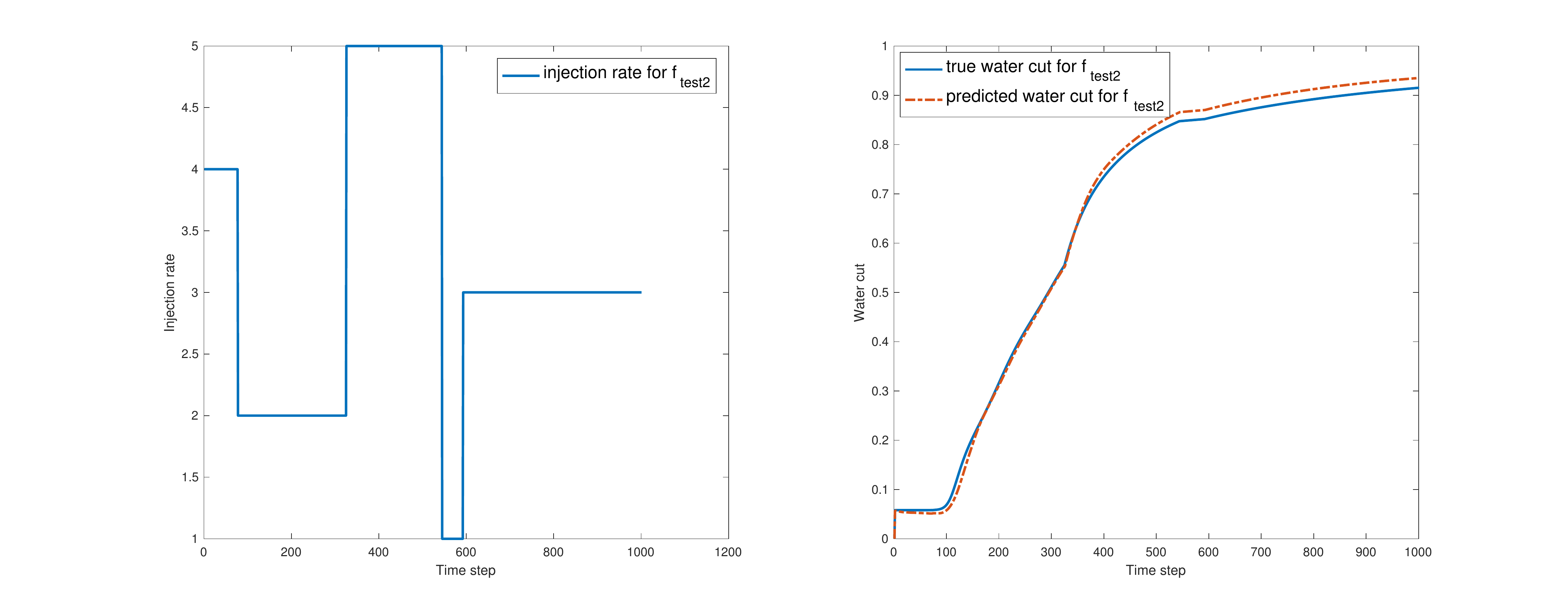}\\
\includegraphics[scale=0.35]{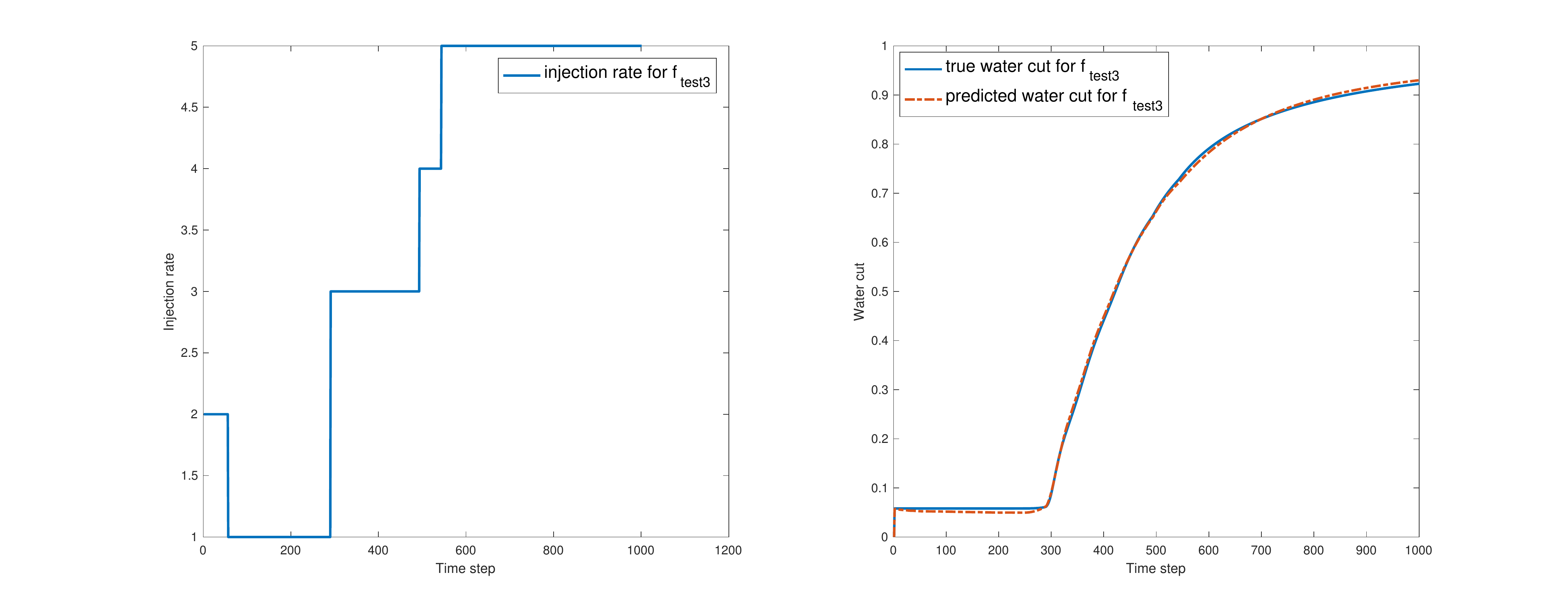}
\caption{ For three different source terms, the comparison of water cuts obtained from the true saturation and the predicted saturation. }
\label{fig:water_cut_tp}
\end{figure}

A visualization of the saturation solutions at some predicted time steps is presented in Figure \ref{fig:sat_nonlinear}, where the source term is the one shown in the left side of Figure \ref{fig:rand_ftest}. We observe that the predicted solutions show a good match towards the true solution. This shows the accuracy of our trained network. In the end, we remark that $999$ steps prediction takes only $2.25$s using the trained networks while solving the coupled system $999$ times takes about $94$s. 
In the end, for three different source terms (shown in left side of Figure \ref{fig:water_cut_tp}), we repeated the prediction as before. Then we also compute their corresponding water cuts based on the true saturation and the predicted saturation. The comparison is presented in Figure \ref{fig:water_cut_tp}. We observe very good match of the predicted and true results.

\begin{table}[!hbt]
\centering
  \begin{tabular}{|c  |c  |}
  \hline
Numbers of predicted time steps& $\norm{S_{\text{pred}} - S_{\text{true}} }_{L^2}$ (\%)  \\  \hline
200  &5.37\\  \hline
400   &4.45	\\  \hline
600   &4.55	\\  \hline
800   &5.42	\\  \hline
1000   &6.88 \\  \hline
  \end{tabular}
\caption{Two phase flow: Given an initial state, apply trained neural networks iteratively. Errors between true and predicted saturation at different numbers of time steps.}\label{tab:sat_coupled_2ph}
\end{table}

One last thing to mention is about the computational time. We list the information in Table \ref{tab:compare_time}.

\begin{table}[!hbt]
\centering
  \begin{tabular}{|c  |c  |}
  \hline
Using Algorithm \ref{alg:sequential} & Using Algorithm \ref{alg:nn_2ph_alg}   \\  \hline
  $94.36$ seconds &$2.25$ seconds\\  \hline
  \end{tabular}
\caption{Comparison of computational time using Algorithm \ref{alg:nn_2ph_alg} and Algorithm \ref{alg:sequential}.}\label{tab:compare_time}
\end{table}

\section{Conclusion} \label{sec:conclusion}

In this work, we propose efficient deep neural networks as surrogate models to approximate flow and transport equations. Both single phase and two-phase flow problems are considered. The designed deep neural networks honor the sparsity structures of the underlying discrete systems, thus having a significant reduction in the number of trainable parameters compared with fully connected neural networks. To be specific, for the flow problem, we apply deep learning to approximate the map from the source terms (in the single-phase flow case)/ relative permeability fields (in the two-phase flow case) to the velocity solution. The networks are constructed with convolutional and locally connected layers to perform model reductions and are equipped with a custom loss function to impose local mass conservation constraints. The custom loss function helps to maintain the physical property of the velocity solution. Once the velocity fields are obtained, they will be inputs to the saturation equation and drive the transport process. Incorporating the learned velocity fields, as well as the upstream scheme of discrete saturation equation, a residual type of network, is introduced to approximate the dynamics of the saturation. That is, we design custom sparsely connected layers which take into account the inherent sparse interaction between the input (saturation at a previous time instant, as well the velocity fields) and output (saturation at the next time instant). Once the feed-forward map between the solution at two consecutive time steps is trained, the approximated map can be used iteratively many times to predict solutions in the long future. An efficient algorithm is proposed to solve the coupled flow and saturation system using the trained neural networks. Considering the saturation map in two-phase flow is nonlinear, the constructed neural networks show great improvement in computational efficiency. The predicted solutions in our numerical examples also present excellent accuracy. Future work involves extending the method in conjunction with multiscale model reduction algorithms and developing more efficient deep learning tools to solve larger coupled flow and transport systems.

\section*{Acknowledgements}
We gratefully acknowledge the support from National Science Foundation (DMS-1555072, DMS-1736364 and DMS-1821233). We also gratefully acknowledge the support of NVIDIA Corporation with the donation of the Titan Xp GPU used for this research.

\bibliographystyle{siam} 
\bibliography{references}

\end{document}